\documentclass[12pt,leqno,twoside]{amsart}
\usepackage{verbatim}
\newcommand{\picspace}{\vspace{3mm}}

\newcommand{\rz}{\mathbb{R}}
\newcommand{\R}{\mathbb{R}}
\newcommand{\gz}{\mathbb{Z}}
\newcommand{\nz}{\mathbb{N}}
\newcommand{\sd}{\mathbb{S}^3}
\newcommand{\sz}{\mathbb{S}^2}
\newcommand{\s}{\mathbb{S}}

\newcommand{\be}{\begin{eqnarray*}}  % Abk"urzungen f"ur Formelzeilen
\newcommand{\ee}{\end{eqnarray*}}
\newcommand{\ben}{\begin{eqnarray}}
\newcommand{\een}{\end{eqnarray}}

\newtheorem{theorem}{Theorem} % [section]
\newtheorem{lemma}[theorem]{Lemma}
\newtheorem{proposition}[theorem]{Proposition}
\newtheorem{corollary}[theorem]{Corollary}

\theoremstyle{definition}

\theoremstyle{remark}

\newcommand{\ia}{{\rm({\it i\/}\rm)}}
\newcommand{\ii}{{\rm({\it ii\/}\rm)}}
\newcommand{\iii}{{\rm({\it iii\/}\rm)}}
\newcommand{\iv}{{\rm({\it iv\/}\rm)}}

\newcommand{\gd}{^\circ}

\newcommand{\cmc}{\textsc{cmc}~}
\newcommand{\acmc}{a \textsc{cmc}~}

\newcommand{\D}{\Omega}

\newcommand{\G}{\mathcal{G}}
\newcommand{\M}{\mathcal{M}}

\newcommand{\Z}{\mathcal{D}}

\renewcommand{\phi}{\varphi}
\renewcommand{\epsilon}{\varepsilon}

\renewcommand{\baselinestretch}{1.12} %Zeilenabstandsvergroesserung

\begin{document}

\thispagestyle{empty}

\author[Grosse-Brauckmann]{Karsten Gro\ss e-Brauckmann}
\thanks{Supported by SFB 256 at Universit\"at Bonn, and
  NSF grant DMS 94-04278 at UMassAmherst}
\address{Universit\"at Bonn\\ Mathematisches Institut\\ Beringstr.\ 1\\
  53115 Bonn, Germany}
\email{kgb@math.uni-bonn.de}

\author[Kusner]{Robert B. Kusner}
\address{Mathematics Department, University of Massachusetts,
  Amherst MA 01003, USA}
\email{kusner@gang.umass.edu}

\title[Symmetric and embedded CMC Surfaces with Few Ends]{Moduli Spaces of 
  Embedded Constant Mean Curvature Surfaces with Few Ends and Special Symmetry}
\date{October 1996}

\begin{abstract}
  We give necessary conditions on complete embedded \cmc surfaces with
  three or four ends subject to reflection symmetries.  The respective
  submoduli spaces are two-dimensional varieties in the moduli spaces
  of general \cmc surfaces.  We characterize fundamental domains of
  our \cmc surfaces by associated great circle polygons in the
  three-sphere.
\end{abstract}

\maketitle

%%%%%%%%%%%%%%%%%%%%%%%%%%%%%%%%%%%%%%%%%%%%%%%%%%%%%%%%%%%%%%%%%%%%%%%%%%%%%%%
%\section{Introduction}

We are interested in explicitly parametrizing the moduli space $\M_{g,k}$ of
complete, connected, properly embedded surfaces in~$\rz^3$ with finite genus
$g$ and a finite number~$k$ of labeled ends~$k$ having nonzero constant mean
curvature.  By rescaling we may assume this constant is~$1$, the mean
curvature of the unit sphere.  Two surfaces in $\rz^3$ are identified as
points in $\M_{g,k}$ if there is a rigid motion of~$\rz^3$ carrying one
surface to the other.  Moreover, we shall include in $\M_{g,k}$ a somewhat
larger class of constant mean curvature (\textsc{cmc}) surfaces, the
\emph{almost} (or \emph{Alexandrov}) \emph{embedded} surfaces, which are
immersed surfaces bounding immersions of handlebodies into~$\rz^3$.  The
space $\M_{g,k}$ is a finite dimensional real analytic variety, and in a
neighborhood of a surface with no $L^2$-Jacobi fields it is a 
($3k-6$)-dimensional manifold~\cite{kmp} for each $k\ge 3$.

A few of the \textsc{cmc} moduli spaces $\M_{g,k}$ are known explicitly: the
only embedded compact \cmc surface is a round sphere~\cite{ale}, so
$\M_{g,0}$ is either a point ($g=0$) or empty ($g>0$); $\M_{g,1}$ is empty,
since there are no one-ended examples~\cite{mee}; two-ended examples are
necessarily the Delaunay \emph{unduloids} \cite{kks}, which are simply-periodic
surfaces of revolution whose minimal radius or \emph{neckradius} $\rho \in
(0,\frac 1 2 ]$ parametrizes $\M_{0,2}$, while $\M_{g,2}$ is empty for $g>0$.
The Kapouleas construction~\cite{kap} shows that $\M_{g,k}$ is not empty for
every $k\ge3$ and every~$g$.  Furthermore, an embedded end of \acmc surface
is asymptotically a Delaunay unduloid~\cite{kks}, and this defines in
particular the neckradii and axes of the ends, which are related via a
balancing formula.  We call the surfaces in $\M_{g,k}$ the
\emph{$k$-unduloids of genus~$g$}, or simply \emph{$k$-unduloids} if their
genus is~0.

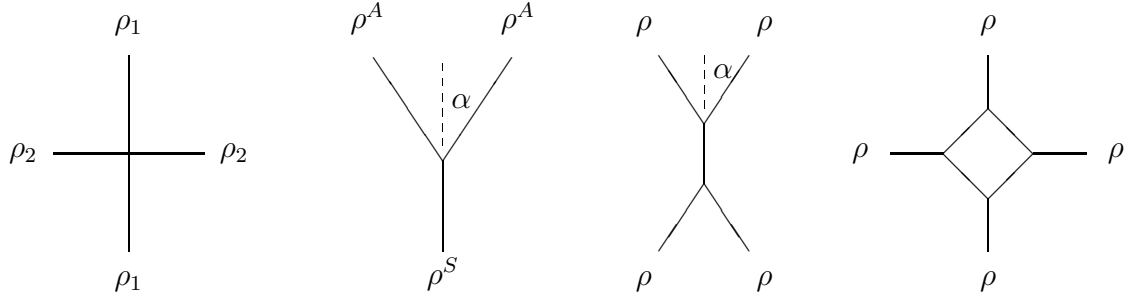
\begin{figure}
  \centerline{
    \unitlength=1mm
    %% rectangular
    \begin{picture}(34.00,40)
      \put(20,8){\line(0,1){26.00}}
      \put(30.00,21){\line(-1,0){20.00}}
      \put(6.00,21){\makebox(0,0)[cc]{$\rho_2$}}
      \put(34.00,21){\makebox(0,0)[cc]{$\rho_2$}}
      \put(20,38){\makebox(0,0)[cc]{$\rho_1$}}
      \put(20,4){\makebox(0,0)[cc]{$\rho_1$}}
    \end{picture}
    \hspace{5mm}
    % isosceles 
    \begin{picture}(30.00,40)
      \put(20.00,20.00){\line(0,-1){12}}
      \put(20.00,20.00){\line(-2,3){9.2}}
      \put(20.00,20.00){\line(2,3){9.2}}
      \put(20.00,5){\makebox(0,0)[cc]{$\rho^S$}}
      \put(30.00,39){\makebox(0,0)[cc]{$\rho^A$}}
      \put(10.00,39){\makebox(0,0)[cc]{$\rho^A$}}
      \linethickness{0.01mm}
      \multiput(20,22)(0,2){6}{\line(0,1){1}}
      \put(22.5,28){\makebox(0,0)[cc]{$\alpha$}}
    \end{picture}
    \hspace{2mm}
    % rhombic
    \begin{picture}(30.00,40)
      \put(20.00,25){\line(0,-1){8}}
      \put(20.00,25){\line(-2,3){6}}
      \put(20.00,25){\line(2,3){6}}
      \put(20.00,17){\line(-2,-3){6}}
      \put(20.00,17){\line(2,-3){6}}
      \put(28,38){\makebox(0,0)[cc]{$\rho$}}
      \put(12,38){\makebox(0,0)[cc]{$\rho$}}
      \put(28,4){\makebox(0,0)[cc]{$\rho$}}
      \put(12,4){\makebox(0,0)[cc]{$\rho$}}
      \linethickness{0.01mm}
      \multiput(20,27)(0,2){4}{\line(0,1){1}}
      \put(22.5,32){\makebox(0,0)[cc]{$\alpha$}}
    \end{picture}
    \hspace{5mm}
    %% rectangular
    \begin{picture}(34.00,40)
      \put(20,8){\line(0,1){7}}
      \put(20,27){\line(0,1){7}}
      \put(14,21){\line(-1,0){7}}
      \put(26,21){\line(1,0){7}}
      \put(20,15){\line(1,1){6}}
      \put(20,15){\line(-1,1){6}}
      \put(20,27){\line(-1,-1){6}}
      \put(20,27){\line(1,-1){6}}
      \put(3,21){\makebox(0,0)[cc]{$\rho$}}
      \put(37,21){\makebox(0,0)[cc]{$\rho$}}
      \put(20,38){\makebox(0,0)[cc]{$\rho$}}
      \put(20,4){\makebox(0,0)[cc]{$\rho$}}
    \end{picture}
    \hspace{5mm}
  }
  \caption{\label{fische}
    Axes and neckradii of rectangular, isosceles, and rhombic surfaces, as 
    well as for a dihedrally symmetric $4$-unduloid of genus~1. 
  }
\end{figure}

In the present paper we focus on \cmc surfaces with few ends and special 
symmetries:  these give two-dimensional submoduli spaces of $\M_{g,3}$ and
$\M_{g,4}$.  The \emph{triunduloids} of genus~$g$ comprising $\M_{g,3}$ are
special in that \emph{a priori} each has a plane of reflection
symmetry~\cite{kks}, which we will think of as a \emph{horizontal} plane.
Their moduli space is 3-dimensional at non-degenerate points~\cite{kmp}.
Here we study the \textsf{Y}-shaped {\em isosceles} triunduloids in
$\M_{0,3}$ which have an additional symmetry in the form of a
\emph{vertical} plane of reflection as depicted in Figure~\ref{fische}(b).
(We will deal with the general triunduloid case elsewhere~\cite{gks}.)  We
also study a submoduli space of the 6-dimensional space $\M_{0,4}$ of
4-unduloids, namely the {\em rectangular} 4-unduloids with the two extra
vertical symmetry planes depicted in Figure~\ref{fische}(a).

%----------------------------------------------------------------------------
%---------------------------------------------------------------------------

One motivation for the present work is the construction of 
a one-parameter family of dihedrally symmetric $k$-unduloids in~\cite{kgb};
these surfaces have $k$~asymptotic axes arranged in a plane with equal
angles~$2\pi/k$.  More precisely, two dihedrally symmetric $k$-unduloids
were shown to exist for neckradii in the interval $(0, \frac{1}{k})$ and one
surface at the right endpoint.  This suggested there
might be some general constraints limiting the allowed neckradius of the ends,
and restricting the number of \cmc surfaces within these limits.  

Indeed, our main result here is a bound on the neckradii for the
two symmetry types (Theorems~\ref{thrsum} and~\ref{thisum}).  These bounds
follow from an analysis of the fundamental domains and their associated
boundary contours.  Using Lawson's theorem \cite{law}~\cite{kar}, if there
is \acmc surface with the assumed reflection symmetries, then there is an
associated minimal surface in the three-sphere~$\sd$ whose boundary consists
of great circle arcs.  The boundary polygons exist only with certain
lengths, and this restricts the range of neckradii the Delaunay ends can
have.  We need to assume enough symmetry to guarantee that a fundamental
domain is simply-connected and Lawson's theorem is applicable.
% but also to have the number of boundary arcs low enough to handle the
% spherical trigonometry explicitly.
%kgb now we know that the period problem is the most important obstacle, 
%    --> [BKS]

We also prove that for given asymptotic axes there are exactly two associated
spherical boundary contours for each neckradius below some maximal value, and
exactly one at the maximum.  We believe (but do not prove) that each contour
bounds a unique minimal surface in~$\s^3$, so that our characterization
of the family of spherical polygons should in fact be a characterization of
the family of \cmc surfaces for the given symmetry type. 

To deduce the existence of a complete \cmc surface from a fundamental
domain may require solving a period problem.  For rectangular 4-unduloids,
the period problem is solved by symmetry alone, while for isosceles
triunduloids the known asymptotics of each end together with the balancing
formula suffice.  This is no longer true for some closely related symmetry
types of surfaces, like the \emph{rhombic} 4-unduloids (see
Figure~\ref{fische}(c)) and dihedrally symmetric $k$-unduloids of
genus~1 (Figure~\ref{fische}(d)); their period problem can still be solved
%kgb: `still' ok?
experimentally~\cite{gbp}.  The rhombic case is particularly interesting
since it allows for a midsection, on which bubbles can be added or deleted
continuously.  A similar property already holds for the isosceles family
(Subsection~\ref{sssmar}):  as one follows a loop in the isosceles
moduli space winding around the point corresponding to a surface with a
cylindrical ``stem'' end, bubbles are generated (or deleted) on the stem.
These observations suggest that the moduli spaces $\M_{0,k}$ may be
connected.

% Finally, we should mention that our approach here depends on the presence of
% reflection symmetries.  This allows us to associate an explicit geometric
% object --- a polygon in $\sd$ --- to a symmetric \cmc surface.  For
% triunduloids, it happens that there is always at least one such reflection
% symmetry, and thus there is hope for extending the methods of this paper to
% study the moduli spaces $\M_{g,3}$ completely --- our work in progress
% \cite{gks} with John Sullivan seems to bear this out in the case of genus~0
% triunduloids.  Of course, almost all \cmc surfaces with $k\ge 4$ ends have
% no symmetries whatsoever, and other approaches will be needed to study the
% \cmc moduli spaces $\M_{g,k}$ in general.  One of us has developed a general
% (though less explicit) approach which associates another geometric object
% --- a \emph{balanced diagram} --- to any \cmc surface; this was explored in
% a very preliminary way in \cite{rob}, and more thoroughly in \cite{kus},
% though it remains a puzzle how to tie the spherical polygon and balanced
% diagram approaches together.

%%%%%%%%%%%%%%%%%%%%%%%%%%%%%%%%%%%%%%%%%%%%%%%%%%%%%%%%%%%%%%%%%%%%%%%%%%
\section{Preliminaries}

%%%%%%%%%%%%%%%%%%%%%%%%%%%%%%%%%%%%%
\subsection{Lawson's theorem}

Our main tool is this consequence of Bonnet's fundamental theorem, as
formulated by Lawson \cite[p.364]{law}:
\begin{theorem}
  Let $M$ be a simply connected immersed minimal surface in $\sd$.  Then
  there exists an isometric immersed \cmc~surface $\tilde M\subset \rz^3$,
  and vice versa. Furthermore, a planar arc of (Schwarz) reflection
  in~$\tilde M$  corresponds to a great circle arc in~$M$.
\end{theorem}
\noindent
We call $M$ and~$\tilde M$ \emph{conjugate cousin surfaces} or simply
\emph{associated surfaces}.  

We will consider fundamental domains of \cmc surfaces with respect to a
group of reflections, whose entire boundary consists in planar arcs of
reflection symmetry (geodesic curvature arcs).
\begin{corollary}
  A simply connected \cmc surface $\tilde M$ bounded by geodesic curvature arcs
  is associated to a spherical minimal surface~$M$ bounded by great circle
  arcs.
\end{corollary}

The fundamental domains of the finite topology \cmc surfaces we study are
not compact, and some bounding geodesic curvature arcs are infinitely long.
We call such an arc a \emph{ray} or a \emph{line} if it extends in one or
both directions to infinity, respectively.  We use the same terminology for
the boundary of the associated minimal domains in $\s^3$:  each arc of the
boundary is a parameterized great circle arc, which in the case of rays and
lines covers the great circle infinitely many times.

%%%%%%%%%%%%%%%%%%%%%%%%%%%%%%%%%%%%%
\subsection{Karcher's Hopf fields}

A~\emph{Hopf field} is defined by the oriented unit tangent vectors to the
great circle fibres of a fixed, oriented Hopf fibration of~$\sd$.  Each
oriented great circle is tangent to a unique Hopf field.  Since there is
an~$\sz$ worth of fibres in a Hopf fibration, and since the space of 
oriented great circles in~$\sd$ is $\sz\times \sz$, there is an~$\sz$ worth
of such fibrations.  Each fibration can be identified with a tangent
direction at some point~$p$ in~$\sd$.  In the following $A$, $B$,~$C$ will
always denote a fixed positively oriented orthonormal basis of the tangent
space $T_p\sd=\R^3$.  With these choices, a Hopf field is determined by a
linear combination $a A+b B +c C$ with $a^2+b^2+c^2=1$.

The spherical boundary contours associated to a \cmc fundamental domain
bounded by geodesic curvature arcs have a property independent of the length
of the arcs:  each arc of the contour determines a Hopf field.  More
precisely the Hopf fields are unique with the following properties:
\ia~fields determined by consecutive great circle arcs make an angle equal
to the angle enclosed by the geodesic curvature arcs at a vertex (this is
the dihedral angle of the two symmetry planes containing the arcs);
\ii\ when measured in terms of Hopf fields, the tangent plane along the
boundary of the spherical minimal surface rotates from one vertex to the
next as much as the tangent plane rotates between two vertices on the \cmc
domain.

Hopf fields were added to Lawson's conjugate surface construction 
by Karcher. See~\cite{kar} or~\cite{kgb} for proofs and details.

%%%%%%%%%%%%%%%%%%%%%%%%%%
\subsection{Ends}

Suppose we divide a \cmc surface by its reflection symmetry group.  If a
symmetry plane contains the axis of an end, we call the portion to one side
of this plane a \emph{half end}.  If another symmetry plane contains the
axis and is orthogonal to the first plane, then the portion contained in the
wedge between these planes is a \emph{quarter end}.

We want to characterize the pair of rays on the boundary of quarter or half
ends.  To do this we look at the case of a quarter or half Delaunay surface
first.
A quarter end of a Delaunay unduloid with neckradius~$\rho$ has an associated
boundary polygon with two great circle rays; the associated spherical minimal
surface domain, which is part of a \emph{spherical helicoid}, contains
perpendicular great circle arcs of length
\ben \label{quarterend}
  \frac{\pi\rho}2 \quad \text{and} \quad \frac{\pi(1-\rho)}2.
\een
Similarly a half end has perpendiculars of length
\ben \label{halfen}
  \pi\rho  \quad \text{and} \quad {\pi(1-\rho)}.
\een
%kgb If you meant this place with . --> ; I don't think it's good.
Since associated surfaces are isometric, this follows from the fact that
there are symmetry circles of radius~$\rho$ and~$1-\rho$ running around a
neck or bulge of the unduloid, respectively.  These circles alternate at 
distance $\pi/2$ for $\rho\not=1/2$; in case~$\rho=1/2$ the surface
associated to the cylinder is foliated by these perpendiculars (it is contained
in a Clifford torus).

An end of an almost embedded \cmc surface with finite topology is
asymptotically a Delaunay unduloid~\cite{kks}, so its neckradius and
axis are defined.  It turns out that the great circle rays bounding 
quarter and half ends coincide with the rays on the boundary associated to
the limiting quarter and half Delaunay ends.
\begin{lemma} \label{leends}
  Let $\tilde M$ be a \cmc fundamental domain which contains a quarter end of
  neckradius~$\rho$.  The associated minimal domain~$M$ in~$\sd$ is then 
  bounded 
  by two great circle rays.  The asymptotic limit of the surface contains 
  perpendicular great circle arcs of lengths given as in~\eqref{quarterend}.
  Similarly, a half end of neckradius~$\rho$ is bounded by two great circle
  rays, and its asymptotic limit contains perpendiculars of lengths
  given as in~\eqref{halfen}.
%
% For $\rho\not=1/2$ there are exactly four such perpendiculars 
% in\/~$\sd$, which are unique up to the 
% antipodal map; in case $\rho=1/2$ the perpendiculars foliate the end, with
% a leaf passing through each point of the boundary rays, and the spherical
% helicoid domain lies on a Clifford torus.
%
\end{lemma}
\begin{proof}
  By Lawson's theorem~$M$ has great circle boundary rays.  Each embedded end
  of~$\tilde M$ is exponentially asymptotic (in the $C^\infty$ topology -- see
  Subsection~\ref{sscont} below) to a Delaunay unduloid~\cite{kks}.  Thus
  $M$ is asymptotic to a spherical helicoid.  A great circle ray asymptotic
  to a great circle must actually cover the great circle.  Hence the
  boundary rays of~$M$ are the great circle rays bounding the associated
  limiting spherical helicoid.
\end{proof}

%%%%%%%%%%%%%%%%%%%%%%%%%%%%%%%%%%%%%%%%%
\subsection{Continuity of families}\label{sscont}

We define a topology on the set of great circle $k$-gons with fixed Hopf fields 
by taking the set of arclengths in~$\rz^k_+$ as continuous coordinates.  We will
also need a topology on families of polygonal contours containing a consecutive
pair of rays.  The pair will always have a shortest perpendicular tangent to
a given Hopf direction.  In this case we can truncate the rays at the
endpoints of this perpendicular, and include their truncated lengths
(modulo~$2\pi$), as well as the length of this perpendicular into the
coordinates.

For families of \cmc surfaces there are various choices for a topology.  If we 
are given a \cmc surface we can write sufficiently $C^0$-close \cmc surfaces as
graphs over compact subsets and use the $C^1$-norm of these graphs.  By
elliptic regularity theory, the resulting topology is equivalent to the
$C^\infty$-topology.  Thus from Lawson's theorem we obtain the following
fact.

\begin{lemma} \label{lecont}
  Consider a family of \cmc surfaces invariant under a fixed group of
  reflections, such that their fundamental domains are simply connected and
  bounded by geodesic curvature arcs.  The map, which assigns to these
  surfaces the great circle polygons bounding the associated minimal
  domains, is continuous.
%kgb  I tried to ungermanize
\end{lemma}

%%%%%%%%%%%%%%%%%%%%%%%%%%%%%%%%%%%%%%%%%
\subsection{Balancing and Kapouleas' existence result}

Let us associate a~{\em force} vector
\ben \label{glforc}
  f := 2\pi\rho(1-\rho) a \in\rz^3
\een
to a Delaunay unduloid with neckradius~$\rho$, whose axis points in the 
direction of the unit vector~$a$. Similarly, for each end of a surface 
in~$\M_{g,k}$, the asymptotic Delaunay limit defines neckradius and axes, 
and hence forces $f_1,\ldots,f_k$.  The \emph{balancing formula}
\ben \label{glbala}
  \sum_{i=1}^k f_i = 0,
\een
is a necessary condition, which in fact holds in much greater 
generality~\cite{kks}. 

It is a deep fact that balancing is also, in a certain sense, a sufficient
condition for the existence of \cmc surfaces: among Kapouleas'
results~\cite{kap} is that $k$-unduloids exist for a dense set of asymptotic
axis directions whose force vectors $f_1,\ldots,f_k$ are balanced and
small.  This smallness condition means that the neckradii are also small.
Here our aim is to study phenomena which occur for large neckradius.

%%%%%%%%%%%%%%%%%%%%%%%%%%%%%%%%%%%%%%%%%%%%%%%%%%%%%%%%%%%%%%%%%%%%%%%%%%%%%%%%
\section{Reflection symmetry}

%%%%%%%%%%%%%%%%%%%%%%%%%%%%%%%%%%%%%%%%%
\subsection{Alexandrov reflection}

% all surfaces are connected!

A.~D.~Alexandrov~\cite{ale} invented a reflection method which was first applied
to show that the round sphere is the only \emph{compact} (almost) embedded \cmc
surface.  This reflection principle was generalized to
noncompact \cmc surfaces in~$\R^3$ as follows~\cite{kks}:
\begin{theorem} \label{legrap}
  Let $M$ be an almost embedded, complete \cmc surface of finite topology.
  If~$M$ lies in a planar slab of
  %kgb
  finite thickness, then there is a plane $P$ in this slab which is a symmetry 
  plane for~$M$.  Furthermore, each symmetric open half $M^+\!$ and~$M^-\!$ 
  of~$M$  is locally a graph lying to one side of~$P$.  In particular the normal
  takes values in open hemispheres on each open half of~$M$.  More precisely, 
  viewing~$M$ as the isometric immersion~$I$ of an abstract surface~$\Sigma$, 
  we can take $\Sigma=\bar\Sigma^+\cup\bar\Sigma^-$ with
  $M^{\pm}=I(\Sigma^{\pm})$.  The restriction of~$I$ to~$\Sigma^{\pm}$ can
  then be written in the form~$I=(\iota,\pm u)$ where $u$ is a positive
  function over some surface domain~$U$, and $\iota$ is an isometric
  immersion of~$U$ into~$P$.
  In case~$M$ is embedded, we may take $M=\Sigma$ and $U\subset P$,
  so that each half~$M^{\pm}$ is globally a graph.
\end{theorem}
\noindent
Let us prove that $M^{\pm}$ lies to one side of~$P$;
all other statements are proved in Remark~2.13 of \cite{kks}.  The boundary
$\partial M^{\pm}$ is contained in~$P$, that is, $u=0$ on~$\partial U$.  By
the asymptotics~$u$ actually takes a minimum.  Suppose that this minimum
value is negative.  Since $M^+$ is locally a graph its mean curvature vector
points downward.  Using a plane parallel to~$P$ as a comparison leads to a
contradiction.

Since each end of an almost embedded \cmc surface with finite topology 
is asymptotically Delaunay, the previous theorem applies when only the
configuration of the axes is known.
%kgb
\begin{corollary} \label{cograp} 
  \ia\
  If the axes of the ends of~$M$ are contained in a slab parallel to a 
%kgb  (Assumptions as in thm.)
  plane~$P$ then $P$ is a symmetry plane. \\
  \ii\
  If $M$ has all its axes parallel, then~$M$ is a Delaunay unduloid.
\end{corollary}
\noindent
For part \ii, the above argument now shows that~$M$ must be cylindrically
bounded, that is, within a bounded distance of a line.  This means we
can apply the reflection argument of the theorem in any direction
about this line, and conclude that $M$ is a surface of revolution
about an axis parallel to this line.

%%%%%%%%%%%%%%%%%%%%%%%%%%%%%%%%%%%%%%%%%
\subsection{Schwarz reflection} \label{sshasz}

If we have a \cmc surface with reflectional symmetries then we can recover
it by Schwarz reflection from a fundamental domain bounded by geodesic
curvature arcs.   In the situation of the above theorem, because $M^{\pm}$
is to one side of the horizontal symmetry plane~$P$, the intersection 
$M\cap P$ consists of geodesic curvature arcs.  We are interested in what
happens when further symmetry planes are present.  In particular we consider
a \emph{vertical} symmetry plane~$V$, which is a plane orthogonal to the
\emph{horizontal} symmetry plane~$P$.

The vertical reflection acts isometrically on the abstract surface~$\Sigma$
with a fixed point set~$F$, and again~$I(F)\subset V\cap M$ consists of
geodesic curvature arcs.  In the embedded case $I(F)=V\cap M$.

When $M$ has genus~0 and $k$~ends we can represent~$\Sigma^+$
(or~$\Sigma^-$, or even~$U$) by a closed disk~$D$ with $k$~points removed 
from~$\partial D$ corresponding to the ends of~$M$.  Furthermore, if we take a
conformal representation, then an isometry of~$\Sigma^+$ acts by a hyperbolic
isometry of~$D$.  Thus the fixed point set $F^+=F\cap\Sigma^+$ of a
reflection is a hyperbolic geodesic.  This geodesic joins either a pair of
boundary punctures, or a boundary point and a puncture, or a pair of
boundary points.  The corresponding cases on~$M^+$ are either a geodesic
curvature line, or a ray, or an arc, respectively.

For genus~1 we can replace the disk representing~$M^+$ with an annulus,
whereas for genus $g>1$ we must allow for $h$~handles attached to the
disk as well as $c$~open disks removed from the interior of~$D$, where
$g=c+2h$.

% note that in general, U need not be a planar domain, if M is not embedded -
% in fact, genus(M)=holes(U)+2handles(U), and, for example, we can make an
% almost embedded cmc surface of genus 2 with handles(U)=1 by taking a balanced
% diagram like this (O's indicate nodes which join):
%
%                       |       |                 
%                       |       |
%                       O       O
%                       |\     /|
%                       | \   / |
%                       |  \ /  |
%                       |   \   |
%                       |  / \  |   U can be thought of as a flat surface which
%                       | /   \ |   thickens this (immersed non-) planar graph
%                       |/     \|
%                       O       O
%                        \     /
%                         \   /
%                          \ /
%                           O
%                           |
%                           |
%                           |

We summarize our discussion above as follows.
\begin{lemma} \label{vertlem}
  For almost embedded $k$-unduloids $M$ of genus~$0$, under the assumption
  of (horizontal) coplanar ends, the fixed point set $F^+$ under reflection in a
  vertical mirror plane~$V$ is connected and is either an arc in $M^+\cap V$
  from the horizontal plane~$P$ to itself, a ray from~$P$ to an end, or a line
  from one end to another.
\end{lemma}

%%%%%%%%%%%%%%%%%%%%%%%%%%%%%%%%%%%%%%%%%%%%%%%%%%%%%%%%%%%%%%%%%%%%%%%%%%
\section{The trigonometry of Lawson quadrilaterals}         \label{setrig}
\subsection{Existence} 

%---------------------------------------------------------------------------
% Der Fundamentalbereich fuer Lawsons Flaeche
%
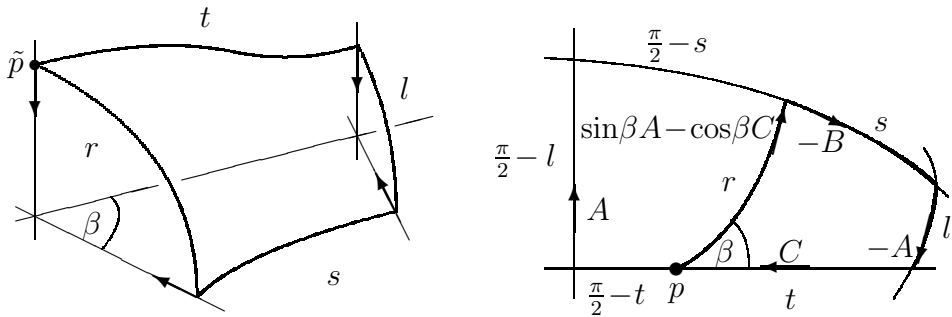
\begin{figure}[bht]
  \centerline{
    \unitlength=0.63mm
    %\begin{picture}(89,62)
    \begin{picture}(89,62)(6,0)
      \thinlines
      \put(10.00,57.00){\line(0,-1){42.00}}
      \put(78.00,32.00){\line(0,1){29}}
      \put(76.00,41.00){\line(1,-2){13.00}}
      \put(6.00,22.00){\line(2,-1){44.00}}
      \put(86.00,39.00){\line(4,1){8.00}}
      \put(6.00,19.00){\line(4,1){30.00}}
      \put(41.00,27.75){\line(4,1){41.00}}
      \put(22.00,18.00){\makebox(0,0)[cc]{$\beta$}}
      %\put(87,47.00){\makebox(0,0)[cl]{$\gt_0,\; -90\gd$}}
      \put(87,47.00){\makebox(0,0)[cl]{$l$}}
      %\put(46,62){\makebox(0,0)[cc]{$\gt_1,\;0\gd$}}
      \put(46,62){\makebox(0,0)[cc]{$t$}}
      %\put(20,39){\makebox(0,0)[cc]{$\gt_2,$}}
      %\put(24,30){\makebox(0,0)[cc]{$-90\gd$}}
      \put(22,34){\makebox(0,0)[cc]{$r$}}
      %\put(73,4){\makebox(0,0)[cc]{$\gt_3,\;90\gd\!-\!\beta$}}
      \put(73,7){\makebox(0,0)[cc]{$s$}}
      \bezier{40}(26,24)(30,20)(24,13) %  um Winkel beta
      \thicklines
      \put(10.00,52.00){\vector(0,-1){11.00}}
      \put(78.00,56.00){\vector(0,-1){10.00}}
      \put(86.00,21.00){\vector(-1,2){4.00}}
      \put(44,3){\vector(-2,1){9.00}}
      \bezier{100}(44,3)(54,13)(86,21)            % untere Kurve, s
      \bezier{70}(10,52)(34,58)(50,55)            % obere Kurve, t, linker Teil
      \bezier{50}(50,55)(65,51.75)(78,56)         % obere Kurve, t, rechter Teil
      \bezier{120}(10,52)(44,35)(44,3)            % linke Kurve, r
      \bezier{60}(78,56)(86,40)(86,21)            % rechte Kurve, l
      \put(10,52){\circle*{2}}
      \put(6,52){\makebox(0,0)[cc]{$\tilde p$}}
    \end{picture}
    \hspace{10mm}
    % Viereck zum Winkel Beta
    \unitlength=0.8mm
    %\raisebox{-12mm}
    {
    %\begin{picture}(82,65) orig
    \begin{picture}(69,43)(15,18)
      \put(15.00,25.00){\line(1,0){65.00}}
      \put(20.00,20.00){\line(0,1){45.00}}
      \bezier{440}(15.00,60.00)(60.00,58.00)(82.00,37.00)
      \bezier{180}(73.00,20.00)(84.00,37.00)(78.00,45.00)
      \put(45.00,27.5){\makebox(0,0)[cc]{$\beta$}}
      \bezier{50}(49.00,25.00)(49.00,31.00)(46.00,33.00)
      \put(17,43){\makebox(0,0)[rc]{$\frac{\pi}2\!-l$}}
      \put(47.00,39){\makebox(0,0)[rc]{$r$}}
      \put(82.00,32){\makebox(0,0)[cc]{$l$}}
      \put(71,49){\makebox(0,0)[cc]{$s$}}
      \put(56,20){\makebox(0,0)[cc]{$t$}}
      \put(27,20){\makebox(0,0)[cc]{$\frac{\pi}2\!-\!t$}}
      \put(37,62){\makebox(0,0)[cc]{$\frac{\pi}2\!-\!s$}}
      \thicklines
      %\put(64,30){\makebox(0,0)[cc]{$\g_1$}}
      \put(59,25.3){\vector(-1,0){8}}
      \put(56,28){\makebox(0,0)[cc]{$C$}}
      \put(20,31){\vector(0,1){8}}
      \put(24,35){\makebox(0,0)[cc]{$A$}}
      %\put(75,35){\makebox(0,0)[cc]{$\g_0$}}
      \put(79.2,34){\vector(-1,-4){2}}
      \put(72.5,28.5){\makebox(0,0)[cc]{$-A$}}
      %\put(57,39){\makebox(0,0)[rc]{$\g_2$}}
      \put(53,44){\vector(1,4){2}}
      \put(53.4,48){\makebox(0,0)[rc]{$\sin\! \beta A\! - \!\cos\! \beta C$}}
      %\put(69,42){\makebox(0,0)[cc]{$\g_3$}}
      \put(58,52){\vector(2,-1){7}}
      \put(61,46){\makebox(0,0)[cc]{$-B$}}
      \put(37.00,25.00){\line(1,0){39.00}}
      \bezier{144}(37.00,25.00)(50.00,32.00)(55.00,53.00)
      \bezier{60}(76.00,25.00)(79.00,30.00)(80.00,39.00)
      \bezier{120}(80.00,39.00)(73.00,46.00)(55.00,53.00)
      \put(37,25){\circle*{2}}
      \put(37,21){\makebox(0,0)[cc]{$p$}}
    \end{picture}
    }
  }
  \picspace
  
  \caption{\label{fidela}
    Quadrilateral bounding a Euclidean \cmc fundamental domain 
    and an isometric spherical minimal domain
    (righthand sketch is meant as a stereographic projection).
  }
\end{figure}
%-----------------------------------------------------------------------------

Our analysis rests on spherical trigonometry, namely the trigonometry of
the boundary polygons of the spherical minimal surfaces associated to the
fundamental domains for the \cmc surfaces with few ends.  The basic polygon
for our analysis is a quadrilateral in $\sd$ with Hopf fields $-A$, $C$, 
$\sin\beta\,A-\cos\beta\,C$, $-B$. Such a quadrilateral has three right 
angles, and one oriented angle $\beta\in[0,2\pi)$
enclosed by the $C$- and $(\sin\beta\,A-\cos\beta\,C)$-arc at some point~$p$.

We will use the letters $l$, $t$, $r$, $s$ to denote the lengths of
a quadrilateral as in Figure~\ref{fidela}.  We call quadrilaterals with the
described Hopf fields, positive lengths, and $0<l\le\pi/4$ \emph{Lawson
quadrilaterals}.  We denote them $\Gamma(l,t,r,s;\beta)$.  They arise from
the truncation of associated quarter ends with an asymptotic perpendicular
of length~$l$ given by Lemma~\ref{leends}; if we choose the shorter
perpendicular then indeed $0<l\le\pi/4$.  The Lawson quadrilaterals are the
associated boundaries of \cmc fundamenal domains coming from doubly periodic
\cmc surfaces provided $\beta=\pi/3$, $\pi/4$, or
$\pi/6$~\cite{law}~\cite{kgb}.

% The two quadrilaterals $\Gamma(\pi/6,\pi/4,\pi/4,\pi/6;\pi/3)$ and
% $\Gamma(\pi/8,\pi/4,\pi/4,\pi/8;\pi/4)$ give rise to doubly periodic \cmc
% surfaces studied by Lawson~\cite{law}.  

In this section we classify all Lawson quadrilaterals.
We first construct two basic families of quadrilaterals with edge lengths
at most $\pi$ and then extend them to obtain all quadrilaterals.
The family \ia\ in the following lemma was already described in Lemma~3.1 
of~\cite{kgb}.
\begin{lemma} \label{letri1}
  For each $\beta\in(0,\pi/2)\cup (\pi/2,\pi)$ there exists a continuous 
  one-parameter family of
  Lawson quadrilaterals $\Gamma(l,t,r,s;\beta)$ as follows:\\
  \ia\ If $0<\beta<\pi/2$ the lengths $l,t,r,s$ range from $0,0,0,0$ past
    $\beta/2,\pi/4,\pi/4,\beta/2$ to $0,\pi/2,\pi/2,\beta$.  Here $r$,
    $s$, $t$ are monotonic, whereas $l$ monotonically increases on the
    first part of the family and monotonically decreases on the second half.
    \\
  \ii\ If $\pi/2<\beta<\pi$ the lengths $l,t,r,s$ range from $0,\pi,0,\pi$
    past $\beta/2,3\pi/4,\pi/4,\pi-\beta/2,$ to $0,\pi/2,\pi/2,\pi-\beta$.
    Again $r$, $s$, $t$ are monotonic, whereas $l$ monotonically increases on
    the
    first part of the family and monotonically decreases on the second half.
\end{lemma}
\begin{proof}
  A Lawson quadrilateral with $0<l,t,s<\pi/2$ and $0<r,\beta<\pi$,
  is uniquely characterized by the four formulas
  \ben
    \label{gleins} \cos s \cos r & = & \cos l \cos t,\\
    \label{glzwei} \sin s \cos r & = & \sin l \sin t,   \\
    \label{gldrei} \cos s\cos l & = & \cos r\cos t + \cos\beta \sin r\sin t,\\
    \label{glvier} \sin s \sin l & = & \cos r \sin t - \cos\beta \sin r\cos t.
  \een
  %
  % The following problem is not addressed:  these equations do not deal
  % with the oriented angle beta.  In fact we have to show directly that
  % a quadrilateral with field  sin beta A - cos beta C  closes for all
  % lengths and beta less than  pi/2  (but does not for -beta) .  
  % To prove this, take a great sphere, containing the  C  and  -A -arc.
  % Provided the length  l  is less than  pi/2  the  -B -arc goes up, say.
  % This also applies to the  sin beta A - cos beta C -arc.  Thus if these
  % two arcs have length at most  pi , then they can match, but the arc
  % with negative beta can not match.
  %                        -B
  %                    __/-->+
  %                 __/      |
  %               _/         |
  %              /           V -A
  %             /            |
  %        -----.-----<------X
  %             A     C     -B
  %
  These four formulas can be obtained using the spherical cosine law, and each 
  is implied by the remaining three, so that one 
  parameter is free.  Indeed in \cite{kgb} it is shown that 
  for $0<l,t,r,s,\beta<\pi/2$ the three equations
  \begin{align} \label{glnull}
    \tan s & =  \tan l \tan t  ,\\
    \label{glacht} 
      \tan 2t& =   \cos\beta \tan 2r ,\\
    \label{glzehn} \quad 
      \cos 2l & %=\cos\beta \sin 2r\sqrt{1+\frac{1}{\cos^2\!\beta\tan^22r}}
      =  \sqrt{\cos^2\!\beta \sin^2 2r + \cos^22r},
  \end{align}
  are equivalent to \eqref{gleins} -- \eqref{glvier} except in
  case $r=\pi/4 \Leftrightarrow t=\pi/4 \Leftrightarrow l=\beta/2
  \Leftrightarrow s=\beta/2$.  For all other
  $0 < r < \pi/2$ we get a length $0 < l \le \beta/2$ by \eqref{glzehn},
  $0 < t < \pi/2$ by \eqref{glacht} and $0 < s < \pi/2$ by \eqref{glnull}.
  This gives a continuous family
  of quadrilaterals for which \eqref{gleins} -- \eqref{glvier} hold,  
  so that part \ia\ is proved.
  % Note that \eqref{glzehn} implies $l \le \beta/2$.

  To prove \ii\ note that 
  $0<l,t,r,s,\beta<\pi/2$ solve \eqref{gleins} -- \eqref{glvier} if and only if
  $\Gamma(l,{\pi-t},r,{\pi-s};\linebreak[0] {\pi-\beta})$ 
  is a Lawson quadrilateral.   To see this we extend 
  Figure~\ref{fidela}(b) to the left up to the next perpendicular of
  length~$l$.  Then the extended part plus the portion to the left of
  the bold quadrilateral form together a large quadrilateral.  Up to the
  motion $(A,B,C)\mapsto(A,-B,-C)$ this quadrilateral has exactly the Hopf
  fields of the bold quadrilateral, but $\pi/2<\beta<\pi$.   Thus it is 
  a Lawson quadrilateral, and the claims for \ii\ follow from~\ia. 
  %
  % The left quadr. has -A  -C  sin beta A - cos beta C  B
  % by the 180 deg. rotation  A,B,C --> A,-B,-C this is the same as
  %                     -A   C  sin beta A + cos beta C -B
  % and the diagonal is  sin (pi-beta) A - cos (pi-beta) C .
  % Thus we arrived as at a quadr. like the original one except that the
  % angle is now  pi-beta .
  %
\end{proof}

The case $\beta=\pi/2$ is somewhat different.  
The equations \eqref{gleins} -- \eqref{glvier} can be checked directly to
confirm the existence of the following Lawson quadrilaterals with 
lengths at most~$\pi$.
\begin{lemma} \label{letri2}
  There exist right-angled Lawson quadrilaterals $\Gamma(l,t,r,s;\pi/2)$ \\
  \ia\ with $0<l\le \pi/4$, $r=\pi/2-l$, and constant lengths $s=t=\pi/2$, \\
  \ii\ with $0<l=r\le \pi/4$, and $s=t=\pi$, and\\
  \iii\ with $l=r=\pi/4$, and $0<s=t$.
\end{lemma} 
These right-angled quadrilaterals bifurcate, for instance 
at $\Gamma(\pi/4,\pi/2,\pi/4,\pi/2;\pi/2)$, which is of type \ia\ and \ii.
As Figure~\ref{fidela} indicates, the right-angled Lawson quadrilaterals are the
associated boundaries of fundamental domains for the Delaunay unduloids.
Therefore these bifurcations can be given the following geometrical 
interpretation.  For the non-cylindrical unduloids the 
planar curvature circles are spaced with distance~$\pi/2$ along the
meridians.  Small and large circles alternate, so that at distance~$\pi/2$ 
there are two different circles (case~\ia\ of the Lemma), and at distance~$\pi$
two equal ones (case~\ii).  On the cylinder, there is a continuous family
of planar symmetry circles having any distance $s=t>0$ (case~\iii).

%%%%%%%%%%%%%%%%%%%%%%%%%%%%%%%%%%%%%
\subsection{Uniqueness} 

We now want to characterize all Lawson quadrilaterals.
In the general case $l,t,r,s\in\rz_+$ 
the `quadratic' equations \eqref{glnull} -- \eqref{glzehn} imply 
\eqref{gleins} -- \eqref{glvier} only up to sign.  Checking cases for the
latter set of equations gives the uniqueness result:
\begin{lemma} \label{leuniq}
  All Lawson quadrilaterals with $0\le \beta<\pi$ are generated from those 
  listed in 
  Lemma~\ref{letri1} and~\ref{letri2} by the following substitutions:\\
  \ia\ adding $2\pi n$ for $n\in\nz$ to  $r$, $s$, or $t$;\\
  \ii\ replacing $(s,t)$ by $(s+\pi,t+\pi)$;\\
  \iii\ replacing $(r,t)$ by $(r+\pi,t+\pi)$;\\
  \iv\ replacing $(r,s)$ by $(r+\pi,s+\pi)$.
\end{lemma}

Whereas there are no Lawson quadrilaterals with $\beta=\pi$, the 
quadrilaterals with $\pi<\beta<2\pi$ are obtained by another substitution from 
Lemma~\ref{leuniq}:
\begin{lemma} \label{lelabe}
  Lawson quadrilaterals with $\pi<\beta<2\pi$ are obtained from 
  Lemma~\ref{letri1} and \ref{letri2} using the following substitutions:\\
  \ia\ replacing $(r,\beta)$ by $(2\pi-r,\beta+\pi)$, or\\
  \ii\ replacing $(r,s,\beta)$ by $(\pi-r,s+\pi,\beta+\pi)$.\\
  All such quadrilaterals are obtained from these by the substitutions
  listed in Lemma~\ref{leuniq}.
\end{lemma}

\section{Rectangular surfaces}

Consider 4-unduloids of arbitrary genus with three orthogonal symmetry
planes.  For combinatorial reasons the axes of the ends must all be
contained in one of the planes, which we call horizontal.  On the other
hand, by balancing and Corollary~\ref{cograp}\ii\ no two axes coincide, and 
there are only two ways for a 4-unduloid to have this symmetry:
\begin{itemize}
  \item The axes form a right-angled cross, so that they are contained
    in the intersection of the horizontal with the vertical symmetry planes
    as in Figure~\ref{fische}(a).
    Opposite ends are congruent and thus there 
    are only two neckradii $\rho_1$ and~$\rho_2$.
    Since the axes are at right angles to each other
    we call such 4-unduloids \emph{rectangular}.
  \item All ends are congruent and their axes form an angle
    $\alpha$ or $\pi/2-\alpha$ with the intersection of the two vertical
    symmetry planes as in Figure~\ref{fische}(c).  There can be a finite length
    midpiece in between the intersection of two pairs of the four axes.  We
    call this case \emph{rhombic}.
\end{itemize}
The two cases have the dihedrally symmetric 4-unduloids~\cite{kgb} in common.

We assume our surfaces have genus~0 for the remainder of this section.  
In our conformal model given in Subsection~\ref{sshasz} we can represent
each half of a four-ended surface by a closed disk with four boundary
punctures.  One vertical reflection can either \ia~fix a pair of opposite
punctures and interchange the other, or \ii~induce a transposition of two
pairs of punctures.  The fixed point set under the vertical reflection joins
the fixed punctures in case~\ia, and by Lemma~\ref{vertlem} there is a
geodesic curvature line on each half of the \cmc surface.  In case~\ii\ the
fixed point set is a geodesic curvature arc which under the horizontal
reflection completes to a closed loop.

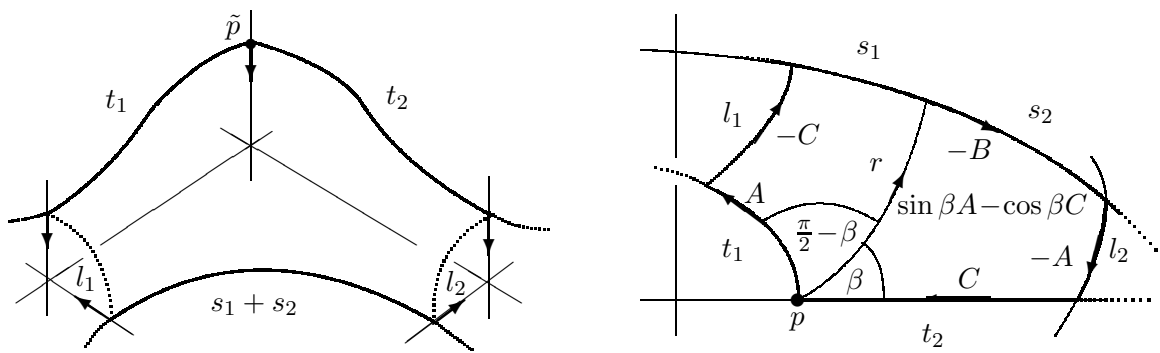
\begin{figure} %[b!]
  \small
  \unitlength=.65mm
  \centerline{
  \begin{picture}(110,69)(12,12)
    % range  12 < x < 122, 12 < y < 81
    \put(61.67,81.33){\line(0,-1){34.67}}
    \put(20.00,46.67){\line(0,-1){27.67}}
    \put(31,33.5){\line(3,2){35}}
    \put(110.33,47.33){\line(0,-1){29.00}}
    \put(95.5,33.5){\line(-5,3){40}}
    \put(62.33,21.33){\makebox(0,0)[cc]{$s_1+s_2$}}
    \put(34.33,63.00){\makebox(0,0)[cc]{$t_1$}}
    \put(91.67,63.67){\makebox(0,0)[cc]{$t_2$}}
    \put(115.33,30.00){\line(-4,-3){22.00}}
    \put(15.67,29.67){\line(3,-2){22.00}}
    \put(28.00,26.67){\makebox(0,0)[cc]{$l_1$}}
    \put(103.67,25.67){\makebox(0,0)[cc]{$l_2$}}
    \put(13.67,21.67){\line(3,2){13.67}}
    \put(117.00,21.67){\line(-5,3){13.67}}
    \thicklines
    \bezier{25}(20.00,39.67)(32.33,34.33)(33.00,18.33)       % l_1-arc
    \bezier{312}(33.00,18.33)(62.67,38.67)(99.00,17.67)      % s arc
    \bezier{25}(110.33,39.33)(97.67,32.67)(99,17.67)         % l_2-arc
    \bezier{16}(99.00,17.67)(102.67,15.00)(106.67,12.00)
    \bezier{12}(33.00,18.33)(30.33,17.00)(27.67,12.67)
    \bezier{140}(110.33,39.67)(92.67,49.33)(84.67,61.67)
    \bezier{112}(84.67,61.67)(80.00,69.67)(62.00,75.00)
    \bezier{108}(62.00,75.00)(52.00,73.33)(41.00,61.00)
    \bezier{124}(41.00,61.00)(32.33,47.67)(19.67,39.67)
    \bezier{16}(110.33,40)(115.00,36.67)(122.33,36.67)
    \bezier{12}(20.00,40.00)(16.00,38.67)(12.00,38.33)
    \put(20.00,40.00){\vector(0,-1){7.33}}
    \put(61.67,74.67){\vector(0,-1){7.67}}                   % tilde p vector
    \put(61.67,74.67){\circle*{2}}                           % tilde p dot
    \put(58,78){\makebox(0,0)[cc]{$\tilde p$}}               % tilde p label
    \put(110.33,39.33){\vector(0,-1){8.00}}
    \put(99.00,17.33){\vector(4,3){6.33}}
    \put(33.00,18.00){\vector(-3,2){7.33}}
  \end{picture}
  \hspace{10mm}
  \unitlength=0.95mm
  % rectangular spherical boundary
  \begin{picture}(69,43)(15,18)
    \put(15.00,25.00){\line(1,0){65.00}}                  % x-line
    \put(20.00,20.00){\line(0,1){21}}                     % y-line, lower part
    \put(20.00,45){\line(0,1){20}}                        % y-line, upper part
    \bezier{440}(15.00,60.00)(60.00,58.00)(82.00,37.00)
    \bezier{180}(73.00,20.00)(84.00,37.00)(78.00,45.00)
    \put(45.00,27.50){\makebox(0,0)[cc]{$\beta$}}
    \bezier{50}(49.00,25.00)(49.00,31.00)(46.00,33.00)
    \put(49,44){\makebox(0,0)[rc]{$r$}}
    \put(82.00,32.00){\makebox(0,0)[cc]{$l_2$}}
    \put(71.00,51){\makebox(0,0)[cc]{$s_2$}}
    \put(56.00,20.00){\makebox(0,0)[cc]{$t_2$}}
    \bezier{144}(37.00,25.00)(50.00,32.00)(55.00,53.00)
    \put(41,34.00){\makebox(0,0)[cc]{$\frac{\pi}2\! -\!\beta$}}
    \bezier{76}(48.00,36.00)(40.00,41.00)(32.00,36.00)
    \put(64.00,25.30){\vector(-1,0){9}}                   % t_2 vector
    \put(61.00,28.00){\makebox(0,0)[cc]{$C$}}
    \thicklines
    \put(79.20,34.00){\vector(-1,-4){1.5}}                % l_2 vector
    \put(72.50,31.00){\makebox(0,0)[cc]{$-A$}}
    \put(50,39){\vector(1,2){2.00}}                       % r vector
    \put(51,38){\makebox(0,0)[cl]{$\sin\beta A\!-\!\cos\beta C$}}
    \put(58.00,52.00){\vector(2,-1){7.00}}                % s vector
    \put(61.00,46.00){\makebox(0,0)[cc]{$-B$}}
    \put(37.00,25.00){\line(1,0){39.00}}                  % t_2 arc
    \bezier{12}(76,25)(79,25)(86,25)                      % t_2 extension
    \bezier{40}(76.00,25.00)(79.00,30.00)(80.00,39.00)    % l_2 arc
    \bezier{140}(80.00,39.00)(65,52)(36,58)               % s arc
    \bezier{12}(80.00,39.00)(82,37.5)(87,32)              % s right extension
    \bezier{12}(36,58)(34,58.5)(29,59)                    % s left extension
    \bezier{96}(37.00,25.00)(37.00,35.00)(24.00,41.00)    % t_1 arc
    \bezier{12}(24.00,41.00)(22,42.5)(17,43.5)            % t_1 extension
    \bezier{80}(24.00,41.00)(36.00,51.00)(36.00,58.00)    % l_1 arc
    \put(31.00,40.00){\makebox(0,0)[cc]{$A$}}
    \put(28.00,32.00){\makebox(0,0)[cc]{$t_1$}}
    \put(28.00,51.00){\makebox(0,0)[cc]{$l_1$}}
    \put(37.00,48.00){\makebox(0,0)[cc]{$-C$}}
    \put(47.00,60.00){\makebox(0,0)[cc]{$s_1$}}
    \put(33,35){\vector(-4,3){7.00}}                      % t_1 vector
    \put(31.00,47.00){\vector(2,3){4}}                    % l_1 vector
    \put(37,25){\circle*{2}}
    \put(37,22){\makebox(0,0)[cc]{$p$}}
  \end{picture}
  }
  \picspace
  
  \caption{ \label{firect}
    Generating \cmc fundamental domain for a rectangular surface,
    and spherical boundary polygon of associated spherical minimal surface.
%   In the asymptotic limit each \cmc end has perpendicular curvature lines
%   of length $l_1$, $l_2$ which are parallel to the two dotted lines.
    The spherical contour is truncated by two asymptotic perpendiculars
    which are associated to the two dotted symmetry lines of length 
    $l_1$, $l_2$ contained in the asymptotic \cmc limit of the end.
  }
\end{figure}

Any other vertical reflection must be of the same type,
unless the surface is dihedrally symmetric.  Thus case~\ia\ corresponds to
the rectangular surfaces, while \ii\ is the rhombic case.  In particular,
a rectangular surface has a fundamental domain bounded by one line and
two rays meeting in a right angle at a point~$\tilde p$ as in
Figure~\ref{firect}(a).  On the other hand, a rhombic surface has fundamental
domain bounded by two arcs and two rays, all meeting at right angles.  The
rectangular case, which we consider here, does not pose a period problem; 
the period problem posed by the rhombic surfaces can be dealt with 
numerically~\cite{gbp}.

%---------------------------------------------------------------------------
% Rectangular 4-und.: truncated CMC fund.domain and its conj. cousin
%
%---------------------------------------------------------------------------

The rectangular fundamental domain contains two quarter ends.
To simplify we use asymptotic perpendiculars provided by Lemma~\ref{leends} to
truncate the spherical contour.   We truncate at the necks, and
by~\eqref{quarterend} the truncated arcs have lengths  
\ben \label{gll1l2}
   0 < l_1, l_2 \le \pi/4. 
\een
The resulting pentagon, with the truncated line of length~$s$ and two truncated
rays of lengths $t_1,t_2$, is
right-angled.  We denote it by $\Gamma(t_1,l_1,s,l_2,t_2)$.

To reduce the spherical trigonometry of these pentagons to the
quadrilaterals of Section~\ref{setrig} we determine the Hopf fields of
the pentagons (Figure~\ref{firect}).
We start at the intersection point~$p$ of the geodesic rays.  The Hopf fields
are then $A$ for one ray, $-C$ for its asymptotic perpendicular, $-B$ for
the line, $-A$ for the other asymptotic perpendicular, and~$C$ for the
returning ray.  We refer to a right-angled
pentagon with these Hopf fields and positive lengths satisfying~\eqref{gll1l2} 
as a {\em truncated rectangular contour}.  The main geometric observation is
a decomposition of these pentagons into quadrilaterals 
as indicated in Figure~\ref{firect}(b):
\begin{lemma} \label{lerepe}
  The set of truncated rectangular contours $\Gamma(t_1,l_1,s,l_2,t_2)$ is
  in 1-1 correspondence to pairs of 
  Lawson quadrilaterals $\Gamma_1(l_1,t_1,r,s_1;\beta)$ and
  $\Gamma_2(l_2,t_2,r,s-s_1;\pi/2-\beta)$ with 
  $0<r,\beta<\pi/2$.
\end{lemma}
\begin{proof}
  If the two quadrilaterals are glued along their arcs of length~$r$
  in such a way that the two complementary angles face each other and are
  contained in the same tangent plane, a right-angled pentagon is formed.
  Moreover, this pentagon has the desired Hopf fields.

  On the other hand, suppose the pentagon is given.  
  We consider a great circle through~$p$ 
  which meets the geodesic containing the opposite $-B$~arc orthogonally in
  a point~$q$.  Extending $s$ and $t_1$, $t_2$ further if necessary, two
  quadrilaterals $\Gamma_1$ and~$\Gamma_2$ result.   
  By orthogonality, the Hopf field of the diagonal is of the form
  $\sin\beta\,A-\cos\beta\,C$ with $\beta\in(0,2\pi)$.  We claim that in
  fact $0<\beta,r<\pi/2$.  

  We use the following fact, which is immediate from Lemmas~\ref{letri1},
  \ref{leuniq}, and~\ref{lelabe}: a Lawson quadrilateral with
  $\beta\in(0,\pi)$ and $\beta\not=\pi/2$ has $0<r\bmod \pi <\pi/2$,
  whereas a quadrilateral with $\beta\in(\pi,2\pi)$ and $\beta\not=3\pi/2$
  has $\pi/2 < r\bmod\pi <\pi$. 

  Suppose that the quadrilateral $\Gamma_2$ has $\beta\in(\pi/2,\pi)$.  This
  means that the quadrilateral~$\Gamma_1$ has angle in $(3\pi/2,2\pi)$.  By
  the preceding fact there is no consistent choice of~$r$ and this case
  is impossible.  The same argument applies to $\beta\in(3\pi/2,2\pi)$.
  The case $\beta\in(\pi,3\pi/2)$ leads to a consistent choice
  $r\in(\pi/2,\pi)$.  However, in this case we can replace the diagonal by an
  arc in the opposite direction; the arc in the opposite direction has
  $\beta\in(0,\pi/2)$ and meets the antipodal point of the former endpoint
  at a length $\pi-r\in(0,\pi/2)$.  Finally, since $l_1,l_2>0$ the 
  angle $\beta$ can
  not be an integer multiple of $\pi/2$, and hence the claim is proved.

  As $\beta\in(0,\pi/2)$, the length of $r$ can be taken to be in
  $(0,\pi/2)$ by Lemma~\ref{letri1}\ia\ and Lemma~\ref{leuniq}. 
  Furthermore, again by Lemma~\ref{leuniq}
  we see that~$q$ is contained in the pentagon and the extension is not
  necessary.
\end{proof}

We now analyze the pentagons via the two quadrilaterals of the lemma.
Note that the family of dihedrally symmetric 4-unduloids
constructed in~\cite{kgb} satisfies $0<\rho\le 1/4$.   These surfaces lead
to symmetric pentagons $\Gamma(l,t,s,t,l)$ with $0<l\le \pi/8$; they
decompose into two equal Lawson quadrilaterals $\Gamma(l,t,r,s/2;\pi/4)$ as
given by Lemma~\ref{letri1}\ia.  For the general case the two 
quadrilaterals are different but $l_1+l_2\le \pi/4$ still holds.

\begin{lemma}
  There exists a continuous two-parameter family of truncated rectangular 
  pentagon contours  
  $${\mathcal G}_{\mathrm{trunc}}=\{ \Gamma(t_1,l_1,s,l_2,t_2) \mid
	 0<l_1,l_2\mbox{\rm \ and }l_1+l_2\le\pi/4,\; 0<t_1,t_2,s<\pi/2 \}. $$
  For each pair $l_1,l_2>0$ with $l_1+l_2<\pi/4$ the family contains
  exactly two distinct contours, while for $l_1+l_2=\pi/4$ there is only one.
\end{lemma}
\begin{proof}
  We determine a quadrilateral 
  $\Gamma_1(l_1,t_1,r,s_1;\beta)$ with angle~$0<\beta<\pi/2$ and
  all edgelengths in~$(0,\pi/2)$ as follows.
  We take $0<l_1<\pi/4$ and $r$ with $l_1 < r < \pi/2-l_1$ as parameters.
  Then \eqref{glzehn} (with $l=l_1$) gives
  % \ben \label{gllofr}
  %   \cos^2\! 2l_1= \cos^2\!\beta\sin^2\! 2r+\cos^2\! 2r ,
  % \een
  \ben \label{glbeta} \cos^2\! \beta =
    \frac{\cos^2\!2l_1-\cos^2\!2r}{\sin^2\! 2r}
  \een
  and we obtain a~$\beta \ge 2l_1$.  Since $l_1<r$ the numerator of
  \eqref{glbeta} is positive, and the fraction is at most 1, so that we can
  choose $\beta<\pi/2$.  Then Lemma~\ref{letri1}\ia\ gives existence of a
  family of quadrilaterals $\Gamma_1(l_1,t_1,r,s;\beta)$ with $0<t_1=t<\pi/2$
  and $0<s<\beta<\pi/2$ as in \eqref{glacht} and~\eqref{glnull}.  These
  quadrilaterals depend continuously on the parameters $r$ and~$l_1$.

  With $r$ and $\beta$ fixed, we now construct another quadrilateral
  $\Gamma_2(l_2,t_2,r,{s-s_1};{\pi/2-\beta})$.
  According to~\eqref{glzehn} and Lemma~\ref{letri1}\ia\ there is such a
  quadrilateral with $0<l_2<\pi/4$ determined by
  \ben \label{glqua2}
    \cos^2\! 2l_2= \cos^2(\pi/2-\beta)\sin^2\! 2r+\cos^2\! 2r ,
  \een
  and $0<t_2<\pi/2$ and $s-s_1<\pi/2-\beta$ given by \eqref{glacht} 
  and~\eqref{glnull}.   Consequently, the quadrilaterals~$\Gamma_2$ form a 
  family which is continuous in $(l_1,r)$.  By Lemma~\ref{lerepe} this 
  gives a truncated rectangular pentagon contour.

  We now want to show that all pairs $l_1$, $l_2$ with $l_1+l_2\le\pi/4$ are
  attained, such that there are two different pentagons for $l_1+l_2<\pi/4$,
  and one if equality holds.  For this we consider the extremal choices of~$r$. 
  If $r\searrow l_1$ then by~\eqref{glbeta} $\beta \nearrow \pi/2$, and thus
  by~\eqref{glqua2} $l_2\searrow 0$.  On the other hand, for $r=\pi/4$,
  \eqref{glbeta} implies $\beta=2l_1$, and from \eqref{glqua2} we conclude
  $2l_2=\pi/2-2l_1$.

  For given $l_1$, it follows from continuity that a choice of~$r$ in the
  lower interval $(l_1,\pi/4]$ gives all $l_2\in(0,\pi/4-l_1]$.  Taking
  $r$-derivatives of \eqref{glbeta} and \eqref{glqua2} it is elementary to
  see that the function $l_2(r)$ is strictly monotonic.  Thus each $l_2$ in
  $(0,\pi/4-l_1]$ is taken exactly once, and in particular $l_1+l_2\le \pi/4$.
  %
  % do it:
  % we look at \eqref{glbeta for r increasing within the interval (l_1,\pi/4):
  % d/dr cos^beta
  %  = [ 2 sin 4r sin^2 2r - (cos^2 2l_1 - cos^2 2r) 2 sin 4r ] / sin^4 2r
  %  = [ 2 sin 4r - cos^2(2l_1) 2 sin 4r ] / ...
  %  = 2 sin^2 2l_1 sin 4r / sin^4 2r = sin^2 2l_1 cos 2r / sin^3 2r  >  0
  % ==> beta(l_1,r) is decreasing in r, beta'<0
  % On the other hand \eqref{glqua2 gives:
  % d/dr(cos^2 2l_2) = cos^2(pi/2-beta)2sin4r + sin(pi-beta)beta'sin^2r - 2sin4r
  %   = - 2sin^2(pi/2-beta)2sin4r + sin(pi-beta)beta'sin^2r  < 0
  % ==> l_2 is increasing

  The upper $r$-interval $[\pi/4,\pi/2-l_1)$ also yields quadrilaterals
  $\Gamma_2$ with $l_2\in(0,{\pi/4-l_1}]$.  This follows with
  the same arguments since
  the limit $r\nearrow \pi/2-l_1$ is analogous to $r\searrow l_1$.
  %
  % if necessary, here is l(r)
  %\ben \label{gllofr}
  %  \cos^2 2l_1= \cos^2\beta\sin^2 2r+\cos^2 2r
  %\een
\end{proof}

Fixing the pentagon at $p$ we extend the lengths $t_1$, $t_2$ and $s$ to
infinity.  We denote by $\mathcal G$ the family of extended contours which 
results this way from ${\mathcal G}_{\mathrm{trunc}}$.  Although there are many
other choices of pentagons for rectangular surfaces besides the family 
${\mathcal G}_{\textrm{trunc}}$, the extended contours are unique.  This
follows from the decomposition of the pentagon into quadrilaterals
(Lemma~\ref{lerepe}), with the Uniqueness~Lemma~\ref{leuniq}.
Furthermore, on the boundary of ${\mathcal G}_{\text{trunc}}$ the lengths
$l_1$ or~$l_2$ vanish.  Thus the family~${\mathcal G}$ is maximal:
\begin{theorem} \label{thrsum}
  The neckradii $\rho_1,\rho_2>0$ of a rectangular $4$-unduloid satisfy
  \ben \label{glrsum}
     \rho_1+\rho_2 \le \frac 1 2 .
  \een
  Furthermore, all rectangular $4$-unduloids have a fundamental domain
  associated to a spherical minimal surface with boundary in the two-dimensional
  family~$\mathcal G$ homeomorphic to an open disk; there are two contours which
%kgb
  satisfy~\eqref{glrsum} with strict inequality, and one for equality.
\end{theorem}

%---------------------------------------------------------------------------
% Rect. moduli space schematically
%
\begin{figure} %[tb]
  \small
  \unitlength=1mm
  \centerline{
  \begin{picture}(47.00,47.00)
    \put(10.00,5){\vector(0,1){34.00}}
    \put(5.00,10.00){\vector(1,0){35.00}}
    %\multiput(35,8)(-4,4){10}{\line(-1,1){3}}
    \bezier{38}(35,10)(22.5,22.5)(10,35)                 % max fam.
    \bezier{80}(17.00,10.00)(16.00,15.67)(10.00,17.00)        % sugg. 48
    \put(16,5){\makebox(0,0)[cc]{\cite{kap}}}            % Nikos
    \put(17,7){\vector(-1,1){4}}
    \put(43,10.00){\makebox(0,0)[cc]{$\rho_1$}}
    \put(10.00,43){\makebox(0,0)[cc]{$\rho_2$}}
    \put(35.00,10.00){\line(0,-1){2.00}}
    \put(35.00,4){\makebox(0,0)[cc]{$\frac 1 2$}}
    \put(10.00,35){\line(-1,0){2.00}}
    \put(5,34.67){\makebox(0,0)[cc]{$\frac 1 2$}}
    \thicklines
    \put(35.00,10.00){\line(-1,1){4}}                     % Jorgen
    \put(35,15){\makebox(0,0)[cc]{\cite{jor}}}
    \put(15,35){\makebox(0,0)[cc]{\cite{jor}}}
    \put(10,35){\line(1,-1){4}}                           % Jorgen
    \put(10.00,10.00){\line(1,1){12.33}}                  % kgb
    \put(22,17){\makebox(0,0)[cc]{\cite{kgb}}}
  \end{picture}
  }
  \vspace{2mm}
  \caption{ \label{firemo}
    The moduli space of rectangular boundary contours can be represented
    by a two-sheeted covering of the triangle defined by 
    $\rho_1,\rho_2 > 0$ and $\rho_1+\rho_2 < 1/2$; the two sheets are
    glued along the dotted line $\rho_1+\rho_2=1/2$.
  }
\end{figure}
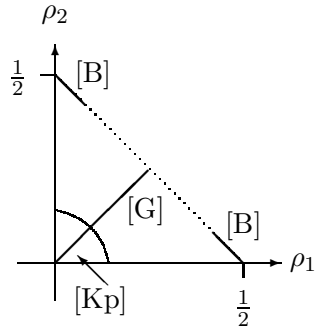
%---------------------------------------------------------------------------

The boundary of the moduli space is given by the contours with
vanishing $\rho_1$ or~$\rho_2$; the degenerate limiting contours bound
surfaces associated to Delaunay surfaces with an orthogonal string of
spheres.  

We believe that our two-parameter family $\mathcal G$ is in 1-1
correspondence to a continuous family of rectangular \cmc surfaces. 
One sheet of the boundary contours with $r>\pi/4$ corresponds to surfaces with
more spherical centers, while on the other sheet with $r<\pi/4$
the centers are smaller.  In particular a rescaled sequence of surfaces 
in the latter sheet with $\rho_1=\lambda\rho_2\to 0$ converges to a minimal
surface with four catenoid ends of alternating logarithmic growth.

To prove existence, one would have to solve Plateau's problem
for the spherical contours in $\mathcal G$ with minimal surfaces.  We believe
these minimal surfaces are unique, at least in the class with almost
embedded associate surfaces.
So far existence is known for the following cases (see Figure~\ref{firemo}):
\begin{itemize}
  \item On an infinite subset with accumulation point the origin
    $\rho_1=\rho_2=0$ by Kapouleas~\cite{kap}.  
%kgb
    This corresponds to the sheet of contours with spherical centers.  It is 
    conceivable that these surfaces form in fact an open continuous family 
    (cf.~Remark~4.6 of~\cite{kap}). 
  \item 
    On the diagonal of $\G$ where the corresponding dihedrally symmetric 
    4-unduloids  attain the neckradii $0<\rho_1=\rho_2\le 1/4$~\cite{kgb}.
  \item 
    Existence of the maximal neckradius family with $\rho_1+\rho_2=1/2$ for
    either $\rho_1$ or $\rho_2$ small was shown by 
    Berglund~\cite{jor}.  Asymptotically, these surfaces are almost 
    cylinders in one direction, and almost spherical unduloids in the
    other.
\end{itemize}
It is known that there are only finitely many components in the moduli space 
of rectangular surfaces with neckradius greater than any fixed
$\epsilon>0$; this follows from the curvature and area bounds of~\cite{koku}. 

%%%%%%%%%%%%%%%%%%%%%%%%%%%%%%%%%%%%%%%%%%%%%
\subsection*{Doubly periodic surfaces of rectangular type}

Using Kapouleas' method~\cite{kap}, doubly periodic surfaces with rectangular
lattice can be found.  Here we want to let the lattice size vary, and
require only that rectangular symmetry is maintained (we still
fix~$H=1$).  Under this assumption Kapouleas has countably many families
${\mathcal F}_{m,n}$, each accumulating at the following degenerate
surfaces: on one set of parallel axes there are $m\ge 0$ spheres in
%kgb
between the junction spheres and likewise $n\ge 0$ spheres in a perpendicular
direction.  
% For these surfaces the two small neckradii $\rho_1$ and~$\rho_2$
% are parameters.  
%kgb
The boundary of the associated minimal fundamental
domains form a family of contours ${\mathcal G}_{m,n}$.  
This family is obtained from ${\mathcal G}_{\text{trunc}}$ as
follows: for a contour in ${\mathcal G}_{\text{trunc}}$ the arc $t_1$ is
extended by $n\pi/2$, $t_2$ by $m\pi/2$, and $s$ by $(m+n)\pi/2$.

Since a continuous family of \cmc surfaces gives rise to a continuous family 
of associated spherical boundary polygons (Lemma~\ref{lecont}) the maximality 
of ${\mathcal G}_{\text{trunc}}$ gives the maximality of ${\mathcal G}_{m,n}$:
\begin{theorem}
%  For each $m,n\ge 0$,
%  Kapouleas' doubly periodic rectangular surfaces ${\mathcal F}_{m,n}$ are
%  contained in a different connected component of the moduli space
%  of doubly periodic surfaces with a set of two orthogonal lattice vectors.
The moduli space
of rectangular doubly periodic surfaces
has connected components indexed by a pair of integers $m,n\ge 0$, each of which contains the corresponding Kapouleas family ${\mathcal F}_{m,n}$.
\end{theorem}
\noindent
Again we believe that the Plateau problem for the contours 
${\mathcal G}_{m,n}$ can be solved with a continuous family of surfaces
extending those of Kapouleas.  The numerical results of~\cite{gbp} make us
doubt that the theorem continues to hold in the class of all (not
necessarily rectangular) doubly periodic surfaces.  Rather, it seems that
only the imposition of a sufficiently strong symmetry group separates the
moduli space into different connected components.

%%%%%%%%%%%%%%%%%%%%%%%%%%%%%%%%%%%%%%%%%%%%%%%%%%%%%%%%%%%%%%%%%%%%%%%%%%%%%%
\section{Isosceles triunduloids}

The class $\M_{g,3}$ of triunduloids of any genus is special in that a
(horizontal) symmetry plane is present~\cite{kks}.  Indeed by balancing
\eqref{glbala} the three force vectors of the ends must be contained in a
plane, which is a symmetry plane by Corollary~\ref{cograp}\ii.  We assume an
additional orthogonal (vertical) symmetry plane and call the triunduloids of
this type \emph{isosceles}.  
% The symmetry implies the axes for the ends of
% an isosceles triunduloid resemble a~\textsf{Y}, and in particular, 
% The axes all meet at a point (this is also true for a general triunduloid 
% and can be proven using torques~\cite{rob}).  
%kgb We can say this in BKS
The intersection of the two symmetry
planes contains the axis of one end, which we call the \emph{stem}. The
other two ends are congruent, and their axes, the \emph{arms}, enclose a
well-defined angle
$$ 
  \alpha\in(0,\pi/2)
$$ 
with the intersection of the symmetry planes (see Figure~\ref{fische}(b)).  
The case $\alpha = 0$ is impossible by Corollary~\ref{cograp}\ia.  Angles
$\pi/2\le\alpha \le \pi$ are excluded by the balancing
formula~\eqref{glbala};  we remark that isosceles triunduloids with
$\pi/2<\alpha<\pi$ exist outside the almost embedded class, for instance
with nodoid ends on either stem or arms~\cite{kap}.
%kgb  (paragr.)

As with the symmetric 4-unduloids in the previous section, we can
characterize a fundamental domain for any isosceles triunduloid of genus~0.
The conformal model in Subsection~\ref{sshasz} lets us represent either half
of a triunduloid as a closed disk with three boundary punctures, and the
vertical reflection must fix one of the punctures and interchange the other
two.  It follows from Lemma~\ref{vertlem} that the fixed point set of the
vertical reflection consists of a geodesic curvature line symmetric
about the horizontal mirror plane and meeting it only at some 
point~$\tilde p$.  Each symmetric half of this line is a ray 
from~$\tilde p$ running out the stem.  Besides one such ray in the vertical
plane, a fundamental domain for any isosceles
triunduloid is bounded by a geodesic curvature ray and line in the 
horizontal plane (see Figure~\ref{fiisos}(a)).
%kgb  (paragr.)

Again the known asymptotics of the ends allows us to truncate the
associated infinite contour with its asymptotic perpendiculars.  
Starting with the point~$p$ associated to~$\tilde p$ (Figure~\ref{fiisos}) the
Hopf fields of the truncated contour can be seen to be $-B$ (ray of half end),
$\cos\alpha\,A-\sin\alpha\,C$ (asymptotic perpendicular for half end), 
$-B$ (line), $-A$ (asymptotic perpendicular of quarter end), $C$ (ray). 
  % proof: along the -B line we must have Hopf field turning alpha
  % denote the field with X, the field capping off the quarter end is -A
  % then for alpha=0:  angle_{-B} (-X,-A) = 0,  i.e. X = A
  %       "        90:                    = 90, i.e. X = -C
  %                     (since -B, C, -A is pos. oriented)
  % Therefore cos alpha A - sin alpha C
We can require
\ben \label{gllrra}
  0 < l \le \frac \pi 4 \qquad \mbox{and} \qquad 0 < r \le \frac\pi 2  ,
\een
if we truncate the contour with the shorter geodesics, i.e.~we cap the ends
at the necks not at the bubbles.  We call a pentagon with the above
Hopf fields which satisfies~\eqref{gllrra} a \emph{truncated isosceles
contour} and denote it with its five lengths $\Gamma(b,r,s-b,l,t)$.

%%%%%%%%%%%%%%%%%%%%%%%%%%%%%%%%%%%%%
\subsection{Trigonometry}

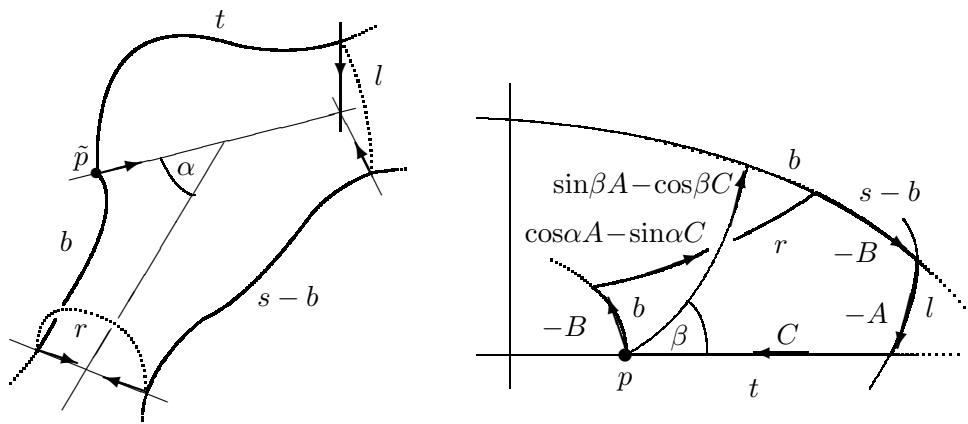
\begin{figure} %[tbh]
  \small
  \unitlength=.5mm
  \centerline{
  \begin{picture}(99,102)(30,22)
  % Euclidean isosceles MC1 patch
  % range  30 < x < 129 , 22 < y < 124
    \put(118.00,94.00){\line(0,1){29.00}}
    \put(116.00,103.00){\line(1,-2){13.00}}
    \put(46.00,81.00){\line(4,1){76.00}}
    % \put(127.00,109.00){\makebox(0,0)[cl]{$\gt_0,\; -90\gd$}}
    \put(127.00,109.00){\makebox(0,0)[cl]{$l$}}
    % \put(86.00,124.00){\makebox(0,0)[cc]{$\gt_1,\;0\gd$}}
    \put(86.00,124.00){\makebox(0,0)[cc]{$t$}}
    \put(77.00,84.00){\makebox(0,0)[cc]{$\alpha$}}
    \qbezier(70.00,87.00)(74.00,78.00)(79,77)
    \put(72.00,22.00){\line(-5,2){42.00}}
    \put(87.00,91.00){\line(-3,-5){28}}
    \put(44,20.00){\line(3,5){13}}
    \thicklines
    \put(118.00,118.00){\vector(0,-1){10.00}}
    \put(126.00,83.00){\vector(-1,2){4.00}}
    \qbezier(90.00,117.00)(52.00,129.00)(53.00,83.00)
    \qbezier(90.00,117.00)(105.00,113.75)(118.00,118.00)
    \qbezier[10](118.00,118.00)(121,119)(127,122)          % r ext of t-arc
    \qbezier[25](118.00,118.00)(126.00,102.00)(126.00,83.00)    % l-arc
    \qbezier[10](126,83)(129,83.5)(135,84)                   % r ext of s-b-arc
    \qbezier(126.00,83.00)(116,79.00)(110,72.00)
    \qbezier(110,72.00)(93,49.00)(81.00,44.00)
    \qbezier(81.00,44.00)(70.00,35.00)(66.00,24)
    \qbezier[8](66,24)(65,21)(65,17)                     % l ext of s-b-arc
    \qbezier[20](66.00,24)(66.00,41.00)(53.00,45.00)        % r-arc, right
    \qbezier[20](53.00,45.00)(37.00,51.00)(37.00,35.5)     % r-arc, left
    \qbezier(37.00,35.5)(41,42.00)(42.00,44)               % b-arc, lower part
    \qbezier[10](37.00,35.5)(35,32.5)(30,27)               % l ext of b-arc
    \qbezier(53.00,83.00)(61.00,74.00)(44.00,48.00)
    \put(45.00,65.00){\makebox(0,0)[cc]{$b$}}
    \put(37.00,36.00){\vector(3,-1){11}}
    \put(66.00,25.00){\vector(-3,1){11}}
    \put(49.00,41.00){\makebox(0,0)[cc]{$r$}}
    \put(104.00,50.00){\makebox(0,0)[cc]{$s-b$}}
    \put(54,83){\vector(4,1){11}}
    \put(53,83){\circle*{3}}
    \put(49,87){\makebox(0,0)[cc]{$\tilde p$}}
  \end{picture}
  \hspace{10mm}
  \unitlength=0.9mm
  \begin{picture}(67,47)(15,18)
  % boundary of spherical isosceles minimal patch
  % range  15 < x < 82  , 18 < y < 65
    \put(15.00,25.00){\line(1,0){65.00}}
    \put(20.00,20.00){\line(0,1){45.00}}
    \bezier{440}(15.00,60.00)(60.00,58.00)(82.00,37.00)
    \bezier{180}(73.00,20.00)(84.00,37.00)(78.00,45.00)
    \put(45.00,27.50){\makebox(0,0)[cc]{$\beta$}}
    \bezier{50}(49.00,25.00)(49.00,31.00)(46.00,33.00)
    \put(82.00,32.00){\makebox(0,0)[cc]{$l$}}
    \put(76,49.00){\makebox(0,0)[cc]{$s-b$}}
    \put(56.00,20.00){\makebox(0,0)[cc]{$t$}}
    \bezier{144}(37.00,25.00)(50.00,32.00)(55.00,53.00)
    \thicklines
    \bezier{12}(76,25)(79,25)(86,25)                      % t_2 extension
    \bezier{12}(80.00,39.00)(82,37.5)(87,32)              % s right extension
    \put(52.5,44.00){\vector(1,4){2.00}}
    \put(49,43){\makebox(0,0)[rc]{$\cos\! \alpha A\! - \!\sin\! \alpha C$}}
    \put(53,50){\makebox(0,0)[rc]{$\sin\! \beta A\! - \!\cos\! \beta C$}}
    \put(64,25.30){\vector(-1,0){8.00}}
    \put(61,28.00){\makebox(0,0)[cc]{$C$}}
    \put(79.20,34.00){\vector(-1,-4){2.00}}
    \put(72.50,31){\makebox(0,0)[cc]{$-A$}}
    \put(37.00,25.00){\line(1,0){39.00}}
    \bezier{60}(76.00,25.00)(79.00,30.00)(80.00,39.00)
    \bezier{48}(37.00,25.00)(37.00,31.00)(32.00,35.00)        % b-arc
    \bezier{12}(32.00,35.00)(30,37)(26,39)        % b-arc left extension
    \bezier{68}(32.00,35.00)(43.00,37.00)(49,40.5)
    \bezier{48}(53,42)(61.00,46.00)(65.00,49.00)   % r-arc upper part
    \put(40.00,37.00){\vector(3,1){8.00}}
    \put(28,29.00){\makebox(0,0)[cc]{$-B$}}
    \bezier{72}(65.00,49.00)(73.00,45.00)(80.00,39.00)  % s-b arc
    \bezier{25}(65,49)(58,52.5)(40,57)           % s-b arc left extension
    \put(62.00,54.00){\makebox(0,0)[cc]{$b$}}
    \put(37,26.00){\vector(-1,3){2.5}}
    \put(39.00,32.00){\makebox(0,0)[cc]{$b$}}
    \put(60,41){\makebox(0,0)[cc]{$r$}}
    \put(70.5,46){\vector(3,-2){8}}
    \put(71,40){\makebox(0,0)[cc]{$-B$}}
    \put(37,25){\circle*{2}}
    \put(37,21){\makebox(0,0)[cc]{$p$}}
  \end{picture}
  }
  \picspace

  \caption{ \label{fiisos}
    Generating \cmc fundamental domain for an isosceles surface with dotted 
    asymptotic curvature lines,
    and spherical boundary polygon of associated spherical minimal surface
    truncated with the asymptotic perpendiculars. 
  }
\end{figure}

Like the rectangular case, the main idea is to reduce the trigonometry
of the pentagons to quadrilaterals.  In this case one of the two quadrilaterals
in the decomposition is what we call a \emph{Clifford rectangle}
$\Gamma_C(b,r)$.  This is a right-angled quadrilateral of edgelength~$b$ and
$r$ which is subset of a Clifford torus.  Opposite arcs are Clifford
parallel, and the Hopf fields can be taken as~$A$, $-B$, $-A$,
$\cos(2b)B-\sin(2b)C$.  Thus in the special case $b=\pi/4$ a Clifford rectangle
is also a rectangular Lawson quadrilateral of the type in
Lemma~\ref{letri2}\iii.  In the general case the Hopf fields
are easily derived from those of a Clifford torus (see~\cite{kgb}).  
\begin{lemma} \label{leistr}
  \ia\ A Lawson quadrilateral $\Gamma_L(l,t,r,s;\pi/2-\alpha+2b)$ and a
  Clifford rectangle $\Gamma_C(b,r)$ can be combined to form a truncated 
  isosceles contour $\Gamma(b,r,s-b,l,t)$ if condition~\eqref{gllrra} and
  the following  hold: $s>b$, and $0<\pi/2-\alpha+2b<\pi$.
  \\
  \ii\ Conversely, a pentagon satisfying \eqref{gllrra} and
  $k\pi<s-b,t<(k+1)\pi$ for sufficiently large $k\in\gz$ can be decomposed
  into a Lawson quadrilateral and a Clifford rectangle as in~\ia.
\end{lemma}
\begin{proof}
  \ia\  If the Clifford rectangle is glued to the Lawson quadrilateral
  such that the $b$-arc of the rectangle is subset of the $s$-arc of
  the quadrilateral, and the $r$-arcs agree (as indicated in 
  Figure~\ref{fiisos}), then a pentagon for isosceles surfaces results.

  \noindent \ii\ 
  We construct a great circle arc through the point~$p$ which meets the
  opposite $-B$ geodesic line orthogonally.
  Since both rays of the half end have the same $-B$~Hopf field,
  they have a continuous set of perpendiculars which are all Clifford 
  parallel to the $r$-arc.  Thus a Clifford parallel to 
  the $r$-arc (in the orientation indicated in Figure~\ref{fiisos})
  at distance~$b$ gives a diagonal with
  Hopf field $\cos(\alpha-2b)\, A-\sin(\alpha-2b)\,C$.
  % original:  $-\sin(\alpha-2b)\, A+\cos(\alpha-2b)\,C$.
  We have to go the same distance $b$ on both $-B$ arcs to obtain the
  Clifford rectangle of the $r$-arc and the diagonal; this means that
  the diagonal only meets an extension of the opposite $-B$ arc.

  Writing the Hopf field of the diagonal in the form
  $\sin(\pi/2-\alpha+2b)\,A -\cos(\pi/2-\alpha+2b)\,C$ we see that
  the other quadrilateral formed has only one non-right angle
  $$
    \beta:=\frac{\pi}2-\alpha+2b
  $$
  at~$p$ and hence is a Lawson quadrilateral.
\end{proof}

%---------------------------------------------------------------------------
% Figure: isosceles fund. domain
%
%---------------------------------------------------------------------------

The problem to determine the possible neckradii amounts to
investigating how the length~$l$ of the asymptotic perpendicular truncating the 
quarter end (stem) relates to the
length~$r$ of the perpendicular truncating the half end (arms).
For the Lawson quadrilateral $\Gamma_L(l,t,r,s;\linebreak[0]{\pi/2-\alpha+2b})$ 
we obtain from \eqref{glzehn} an expression involving the length $b$,
\ben \label{gltrig}
  \cos^2 2l = \sin^2 (\alpha-2b) \sin^2 2r + \cos^2 2r          .
\een
%
% <==>   [ cos^2 2l - cos^2 2r ] / sin^2 2r  =  sin^2 (alpha - 2b)
%      RHS is > 0 for       0 < b < alpha/2  (first branch)
%             <    "  alpha/2 < b < alpha   (second branch)
% ==>                                     /  cos^2 2l - cos^2 2r
%    alpha - 2b = sign(alpha-2b) arcsin  /  ---------------------
%                                      \/         sin^2 2r
%

%%%%%%%%%%%%%%%%%%%%%%%%%%%%%%%%%%%%%
\subsection{Balancing and the period problem}

There is a period problem for isosceles surfaces:  
the associated \cmc domain must have the two rays bounding the half end
in the \emph{same} plane, not just in parallel planes.
Boundary contours which can bound surfaces with vanishing periods form
a codimension one family within all truncated isosceles contours.  
Although we do not give an existence proof, we can select such a family
using two facts: the balancing formula, and the asymptotics of the ends.  
Technically the condition we obtain is a further necessary condition
on the boundary contours.  

The balancing formula~\eqref{glbala} relates the forces $f^A$ of the arms
to the force~$f^S$ of the stem
\be
  |f^S| = 2\cos\alpha \,|f^A|  .
\ee
Thus \eqref{glforc} yields for the asymptotic neckradii $\rho^S$ and~$\rho^A$
\ben \label{glnetr}
    2\pi \rho^S (1-\rho^S) = 2\cos\alpha \, 2\pi \rho^A (1-\rho^A)   .
\een
The arc lengths of the perpendiculars in the truncated contour are related
to the neckradii,
$$2\pi \rho^S = 4 l \qquad \text{and} \qquad 2\pi \rho^A = 2 r,$$
and we obtain
%kgb
% \be
%   4 l \left(1 - \frac{2l}\pi\right) =
%     2\cos\alpha \: 2 r \left(1-\frac r \pi\right)
% \ee
% or, simply,
\ben \label{glbasi}
  l (\pi - 2l) = \cos\alpha \, r (\pi - r).
\een
% i.e.                         l ( pi - 2l)
%                 cos alpha = --------------
%                               r (pi - r)
%
Solving this quadratic equation 
for those~$l$ admissible by~\eqref{gllrra} we arrive at 
the following statement on the pentagon lengths, which, by~\eqref{glnetr} 
is in fact a result on the neckradii.
\begin{lemma}
  Let $\Gamma$ be the contour bounding a minimal surface in~$\s^3$ associated to
  the fundamental domain of a (balanced) isosceles surface.
  If the great circle rays of the half end in $\Gamma$ have 
  a shortest perpendicular with length~$r$,
  then the great circle rays of the quarter end have shortest perpendiculars
  of length 
  \ben \label{glbalr}
    0 < l = \frac \pi 4 - \sqrt{\frac{\pi^2}{16}-\cos\alpha\,\frac{r(\pi-r)}2 }
    \le \frac\pi 4.
  \een
\end{lemma}
% this is defined for
%    pi^2/16-cos alpha r(pi-r)/2 > 0
% <==>  r^2  - pi r  + pi^2/[8 cos alpha] > 0
% ==> this vanishes for  r = pi/2 - sqrt{ pi^2/4 - pi^2/[8 cos alpha] }
%                          = pi/2 (1-sqrt{1-1/[2cosalpha]})
%
In particular, \eqref{glbalr} has only a solution for those $0\le r \le \pi/2$
which are not greater than
$$r_{\max}(\alpha):= \left\{
      \begin{array}{lcl}
	\frac{\pi}2\left(1-\sqrt{1-\frac 1{2\cos\alpha}}\right) \;
	  & \mbox{ for } & 0<\alpha\le\frac{\pi}3   \\
	\frac{\pi}2    & \mbox{ for } & \frac{\pi}3\le \alpha < \frac{\pi}2 \, .
      \end{array}
    \right.
$$

%  \ben \label{glrmax}
%   I(\alpha):= \left\{
%     \begin{array}{lcl}
%       \Big[ 0,\frac{\pi}2\left(1-\sqrt{1-\frac 1{2\cos\alpha}}\right)\; \Big]
%         & \mbox{ for } & 0<\alpha\le\frac{\pi}3   \\
%       \left[0,\frac{\pi}2\right]
%         & \mbox{ for } & \frac{\pi}3\le \alpha < \frac{\pi}2 \, .
%     \end{array}
%   \right.
% \een
%
%
% Here is the inverse function r(l) :
%                                _____________________
%                   pi         / pi^2     l (pi-2l)
%            r(l) = ---  -    /  ----  -  ---------
%                    2      \/    4       cos alpha
%
% The equivalent condition for the radii (check again):
%                  _____________________________________
%          1      / 1      2                   2
%  rho   = - -   /  -  -  --- cos alpha rho (pi -rho )
%     S    2   \/   4      pi              A        A
%

%%%%%%%%%%%%%%%%%%%%%%%%%%%%%%%%%%%%%%
\subsection{Truncated associated boundary contours}

We now study the truncated isosceles contours via the Decomposition 
Lemma~\ref{leistr}.  The first step is to describe a family of quadrilaterals 
satisfying \eqref{gltrig} and~\eqref{glbalr}.  
Eliminating~$l$ from these equations means that $(\alpha,r,b)$
is a zero of the function
\be
  f(\alpha,r,b):=\sin^2(\alpha-2b)\sin^2\!2r+\cos^2\!2r
    -\sin^2\sqrt{\frac{\pi^2}{4}-2\cos\alpha\,r(\pi-r)}.
\ee
Conversely, a zero $(\alpha,r,b)$ of~$f$ and the length $l$ defined
by~\eqref{glbalr} solve \eqref{gltrig}. 
We now want to determine all zeros of~$f$.  By the periodicity of~$f$ in~$b$ 
it will be sufficient to confine $b$ to some interval of length~$\pi$.
Specifically, we set
$$\D:=\left\{(\alpha,r,b) \mid
  0<\alpha<\pi/2,\; 0<r\le r_{\max}(\alpha),\;
  \frac\alpha 2-\frac\pi 4 \le b \le \frac\alpha 2+\frac\pi 4
  \right\} .
$$
%
% We can solve f explicitly for b
%  [ sin^2(2sqrt(...)) - cos^2r ] / sin^2 2r = sin^2 (alpha - 2b)
% this gives                               ________________________________
%         alpha                           /  sin^2(2 sqrt(...)) - cos^2 2r
%    b =  ----- - sign(alpha-2b) arcsin  /   -----------------------------
%           2                          \/               sin^2 2r
%

\begin{lemma}
  On the set $\D$, all zeros of $f$ form a continuous two-parameter family 
  homeomorphic to an open disk
  $$
    \Z := \{ (\alpha,r,b) \in \D \mid
	0 < r \le R(\alpha), b=b_i(\alpha,r)\; \textrm{ for }\; i=1,2 \},
  $$
  where 
  $$ R(\alpha):= \pi\frac{1-\cos\alpha}{2-\cos\alpha}.$$
  % some values: R(0) = 0
  %              R = pi/2 O(alpha^2)    for small argument
  %              R(arccos 2/3) = pi/4
  %              R(pi/2) = pi/2
  %  R'(\alpha) = pi sin alpha / (2-cos alpha)^2
  %               i.e. R'(pi/2) = pi/4
  %  R"(\alpha) = pi 2(cos alpha- sin alpha) - cos^2 alpha / (2-cos alpha)^3
  %               i.e. only convex up to some angle less than 45 degrees
  Here 
  $b_1(\alpha,r)\in(\alpha/2-\pi/4,\alpha/2]$
  and $b_2(\alpha,r)\in [\alpha/2,\alpha/2+\pi/4)$
  are two continuous functions with 
  $b_1(\alpha,R(\alpha))= b_2(\alpha,R(\alpha))=\alpha/2$. 
  % and $\lim_{r\to 0} b_1(\alpha,r)=0$, $\lim_{r\to 0} b_2(\alpha,r)=\alpha$.
\end{lemma}
%---------------------------------------------------------------------------
% Figure: Plot of R(alpha)
%
% known: R(0)=0, R(arccos(2/3)=pi/4, R(\pi/3)=pi/3, R(pi/2)=pi/2
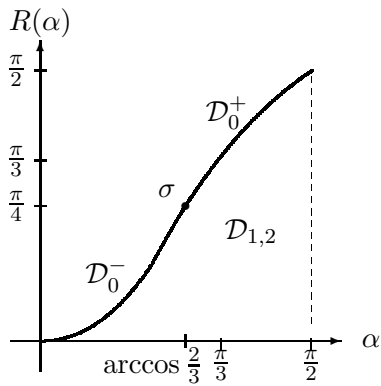
\begin{figure} [bth]
  \small
  \unitlength=.8mm
  \centerline{
  \begin{picture}(60,60)(-5,-7)
    \put(0,-5){\vector(0,1){55}}        % x axis
    \put(-5,0){\vector(1,0){55}}        % y axis
    \put(55,0){\makebox(0,0)[cc]{$\alpha$}}        
    \put(0,53){\makebox(0,0)[cc]{$R(\alpha)$}}      
    \put(45,-1.5){\line(0,1){2}}                     % pi/2 marker  x axis
    %\put(45,-4){\makebox(0,0)[cc]{$\scriptstyle \pi/2$}}      
    \put(45,-4){\makebox(0,0)[cc]{$\frac \pi 2$}}      
    \put(30,-1.5){\line(0,1){2}}                     % pi/3 marker  x axis
    %\put(31,-4){\makebox(0,0)[cc]{$\scriptstyle \pi/3$}}      
    \put(30,-4){\makebox(0,0)[cc]{$\frac\pi 3$}}      
    %\put(22.5,-1){\line(0,1){2}}                     % arccos(2/3)  x axis
    \put(24.1,-1){\line(0,1){2}}                     % arccos(2/3)  x axis
    % arccos(2/3) = 48.2 degrees
   %\put(27,-4){\makebox(0,0)[rc]{$\scriptstyle\arccos{\scriptstyle\frac 2 3}$}}
    \put(27,-4){\makebox(0,0)[rc]{$\arccos{\frac 2 3}$}}
    \put(-1,45){\line(1,0){2}}                     % pi/2 marker  y axis
    \put(-4,45){\makebox(0,0)[cc]{$\frac\pi 2$}}      
    \put(-1,30){\line(1,0){2}}                     % pi/3 marker  y axis
    \put(-4,30){\makebox(0,0)[cc]{$\frac\pi 3$}}      
    \put(-1,22.5){\line(1,0){2}}                     % pi/4 marker x axis
    \put(-4,22.5){\makebox(0,0)[cc]{$\frac\pi 4$}}      
    \put(24.1,22.5){\circle*{1.5}}                      % point (arccos,pi/4)
    \put(21,25){\makebox(0,0)[cc]{$\sigma$}}      
    \thicklines
    \qbezier(0,0)(10,0)(18.5,13)                     % lower arc
    \qbezier(18.5,13)(30,35)(45,45)                     % upper arc
    \put(11,11){\makebox(0,0)[cc]{$\Z^-_0$}}      
    \put(31,38){\makebox(0,0)[cc]{$\Z^+_0$}}      
    \put(35,18){\makebox(0,0)[cc]{$\Z_{1,2}$}}      
    \linethickness{0.01mm}
    \multiput(45,3)(0,2){21}{\line(0,1){1}}        % dashed right line
  \end{picture}
  }
  \vspace{2mm}

  \caption{ \label{fiisor}
     The function $R(\alpha)$ and the projection of $\Z$.  
%    The sets $\Z^\pm$ project onto the area between $x$-axis and graph.
  }
\end{figure}
%---------------------------------------------------------------------------

\begin{proof}
  It is straightforward to show that $f(\alpha,r,b)$ is defined on $\D$,
  i.e.~$ R(\alpha) \le r_{\max}(\alpha)$.
%  Since $R(\alpha)<\pi/2$ this is obvious for $\alpha\ge\pi/3$.  For
% $0<\alpha<\pi/3$ \eqref{glinte} is equivalent to
% $$
%   1-2\;\frac{1-\cos\alpha}{2-\cos\alpha} \ge \sqrt{1-\frac1{2\cos\alpha}}.
% $$
% The left hand side is equal to
% $[2-\cos\alpha-(2-2\cos\alpha)]/[2-\cos\alpha]= \cos\alpha/[2-\cos\alpha]$.
% Since this is positive squaring gives the equivalent inequality
% $$
%   \frac{\cos^2\!\alpha}{\cos^2\alpha-4\cos\alpha+4}
%    \ge \frac {\cos\alpha-\frac 1 2}{\cos\alpha},
% $$
% or $\cos^3\!\alpha \ge %(\cos\alpha-\frac 1 2)(\cos^2\!\alpha-4\cos\alpha+4) =
% \cos^3\!\alpha -\frac 9 2 \cos^2\!\alpha+6\cos\alpha-2$.
% This inequality is valid since
% $$
%   -\frac 9 2 \cos^2\!\alpha + 6\cos\alpha-2
%     = -\left(\frac 3{\sqrt 2}\cos\alpha - \sqrt 2\right)^2\le 0,
% $$
% and holds with equality for $\cos\alpha = 2/3$ only.
%  
  We want to find zeros of $f$ on $\D$.
  We establish $b_1(\alpha,r)$ first for $r$ in the open interval 
  $0<r<R(\alpha)$.
  % which means that $0<\beta<\pi/2$.
  An explicit formula for $b_1(\alpha,r)$ involves many square roots, so 
  it is more 
  straightforward to apply the implicit function theorem.  To do so we claim:\\
  \ia\  $\frac{\partial f}{\partial b} <0 $ on $\D$, \\
  \ii\ $f(\alpha,r,\alpha/2-\pi/4)>0$ for $0<\alpha<\pi/2$ and 
       $0<r\le R(\alpha)$,     and \\
  \iii\  $ f(\alpha,r,\alpha/2)<0$ for $0<\alpha<\pi/2$ and $0<r<R(\alpha)$.

  The derivative
  \ben \label{glbder}
    \frac{\partial f}{\partial b} = -2\sin(2\alpha-4b)\sin^2\!2r
  \een
  is clearly negative on $\D$ which gives claim~\ia.
  To prove \ii\ observe $\sin^2(\alpha-2(\alpha/2-\pi/4))=1$, so we have
  $f(\alpha,r,\alpha/2-\pi/4)
    =\sin^2 \!r+\cos^2 \! r
      -\sin^2\left(2\sqrt{\frac{\pi^2}{16}-\cos\alpha\frac{r(\pi-r)}2}\right)
    = \cos^2(2\sqrt{\ldots})$.
  This is non-negative and vanishes exactly at $r=0$.

  Proving \iii\ amounts to a further elementary calculation:
  $f(\alpha,r,\alpha/2)<0$~is equivalent to
  \be
    \sin^2\left(\frac{\pi}2-2r\right)
    < \sin^2\sqrt{\frac{\pi^2}{4}-2\cos\alpha r(\pi-r)}.
  \ee
  Since $0<\sqrt{\ldots} < \pi/2$ and $-\pi/2<\pi/2-2r<\pi/2$ this is
  equivalent to
  $$
    \left|\frac{\pi}2 - 2r \right|
    < \sqrt{\frac{\pi^2}{4}-2\cos\alpha r(\pi-r)}.
  $$
  Squaring yields
  \be
    \frac{\pi^2}{4}-2r(\pi-2r)
    < \frac{\pi^2}{4}-2\cos\alpha r(\pi-r) ,
  \ee
  which is equivalent to $\pi-2r > \cos\alpha (\pi-r)$,
  $\pi(1-\cos\alpha)>(2-\cos\alpha)r$ or, finally, $0<r<R(\alpha)$.

  For each $\alpha$ the implicit function theorem gives a differentiable 
  function~$b_1(\alpha,r)$ where $0<r<R(\alpha)$.   The function $b_1$ is 
  unique in the subset of~$\D$ with $\pi/4-\alpha/2<b<\alpha/2$ by~\ia.  

  It is a direct consequence of the previous calculation that 
  $(\alpha,R(\alpha),\alpha/2)$ is also a zero of~$f$.  Moreover
  from \eqref{glbder} it follows that this is the only zero in~$\D$ with
  $r=R(\alpha)$ and $b\le\alpha/2$.  Thus we can continuously extend 
  $b_1$ by setting $b_1(\alpha,R(\alpha))=\alpha/2$.
  Furthermore, if $r>R(\alpha)$ it follows from the calculation that 
  $2\sqrt{ \pi^2/{4}-2\cos\alpha r(\pi-r) } < |\pi/2-2r|$;
  thus 
  $f(\alpha,r,b)>\sin^2(\alpha-2b)\sin^2\!2r    % + \cos^2\!2r-\sin^2(\pi/2-2r)
  >0,$
  and there are no zeros of $f$ for $r>R(\alpha)$.

% It remains to prove that % $b_1\ge 0$.
% $b_1\to 0$ for $r\searrow 0$.  This follows from 
% \ben \label{glfior}
%   f(\alpha,b_1,r) = 
%         4r^2 \left(\cos^2\!\alpha-\cos^2(\alpha-2b_1)\right)+ O(r^4)
% \een
% for small $r$.
  %
  % Here is the calculation.  We use  cos^2 x = 1 - x^2 + O(x^4),
  %   and  sqrt(c+x) = sqrt c + x / [2 sqrt c] + O (x^2)
  %
  % f(b,r) = sin^2(alpha-2b)sin^2 2r + cos^2 2r
  %          - sin^2 ( 2sqrt[pi^2/16 - cosalpha r(pi-r)/2 ] )
  %        = sin^2(alpha-2b)(4r^2+O(r^4)) + 1 - 4r^2 + O(r^4)
  %          - sin^2 ( 2 [ pi/4 -  cos alpha r + O(r^2) ] )
  %        = 1 - 4 r^2 + sin^2(alpha-2b) 4 r^2 + O(r^4)
  %          - cos^2 ( 2 r cos alpha + O(r^2) )
  %        = 1 - 4 r^2 cos^2(alpha-2b)
  %          - 1 + 4 r^2 cos^2(alpha)  + O(r^4)
  %

  To obtain $b_2\in [\alpha/2,\pi/4)$ we set
  $b_2(\alpha,r):=\alpha-b_1(\alpha,r)$.
  Then $\alpha-2b_1 = 2b_2-\alpha$ and this leaves $f$ is invariant.
  All properties claimed for $b_2$ follow from $b_1$.
  % Since $0<\beta<\pi/2$ on the first branch
  % we have $\pi/2<\beta<\pi$ on the second branch.
\end{proof} 

To determine the Lawson quadrilaterals in the form for 
Lemma~\ref{leistr}, we have to select values of $l$, $t$, and $s$
satisfying \eqref{gleins} -- \eqref{glvier}.  Whereas there is a continuous
way to choose~$l$, the two lengths $t$ and~$s$ which reflect the
position of the truncation by the asymptotic perpendiculars
must have a discontinuity on~$\Z$.  
In the following lemma we first establish the quadrilaterals on 
$\Z_{i}=\{(\alpha,b,r)\in \Z \mid b =b_i(\alpha,r), r\not=R(\alpha)\}$ 
for $i=1,2$,
%$\Sigma_{1}=\{(\alpha,b,r)\in \Sigma \mid b < \alpha/2\}$,
%$\Sigma_{2}=\{(\alpha,b,r)\in \Sigma \mid b > \alpha/2\}$,
and $\Z_0:= \{(\alpha,\alpha/2,R(\alpha))\}$.
Then we discuss the continuity problem.  It will be useful to decompose
$\Z_0$ further, namely into the two sets 
$\Z_0^{-} := \{(\alpha,b,r)\in\Z_0 \mid \alpha < \arccos(2/3) \}$
and 
$\Z_0^{+} := \{(\alpha,b,r)\in\Z_0 \mid \alpha > \arccos(2/3) \}$,
as well as the point $\sigma:=(\arccos(2/3),\arccos(2/3)/2,\pi/4)$.
We can now state the main technical lemma for the isosceles case.

\begin{lemma}\label{lemain}
  For each point $(\alpha,b,r)\in \Z-\{\sigma\}$ there is a Lawson
  quadrilateral $\Gamma_L(l,t,r,s;\linebreak[0] \pi/2-\alpha+2b)$ satisfying
  \eqref{gltrig}, \eqref{glbalr}, and~\eqref{gllrra}.  According to the
  decomposition of $\Z-\{\sigma\}$ into $\Z_1 \cup \Z_2 \cup \Z_0^{\pm}$ we write
  the families of these quadrilaterals as $\G_1 \cup \G_2 \cup\G_0^{\pm}$.
  Moreover, these quadrilaterals are unique with $0<s,t\le\pi$.  For the
  point~$\sigma$ there is a one parameter family of distinct quadrilaterals
  ${\mathcal G}_0^\sigma$ with $0<s,t\le\pi$.  The length~$l$ of all these
  quadrilaterals induces a continuous function on~$\Z$, which vanishes on
  $\partial \Z$.  The lengths $t,s$ are continuous on $\Z_1 \cup \Z_0^+ \cup
  \Z_2$.  When $(\alpha,r,b)\in\Z_2$ approaches a point in~$\Z_0^-$, we
  have  $\lim s = \lim t = \pi$; whereas for $(\alpha,r,b)\in\Z_1$
  approaching such a point, $\lim s = \lim t = 0$.
%kgb  clearer?
\end{lemma}
\begin{proof}
  A unique value $0<l\le\pi/4$ on~$\Z$ is determined by~\eqref{glbalr}.
  Since $r$ and $\alpha$ are continuous on~$\Z$ so is~$l$.
  On the portion of~$\partial \Z$ with $r=0$ 
  it follows that $l\equiv 0$.
  The remaining boundary of $\Z$ satisfies $\alpha=\pi/2$.  Again,
  it follows directly from~\eqref{glbalr} that $l\searrow 0$ when this
  boundary arc is approached.

  \noindent
  \ia\
  We now want to define $s,t$ on $\Z_1 \cup \Z_2$. 
  For the points in $\Z_1$ with $b_1<\alpha/2$ we
  have $0<\beta<\pi/2$.  In this
  case~\eqref{glacht} and~\eqref{glnull} define $0<s,t<\pi/2$ 
  continuously, and uniquely within $(0,\pi]$.  This gives~$\G_1$.

% When $r\searrow 0$ we proved in the last lemma that $b\searrow 0$
% and thus $\beta\searrow \pi/2-\alpha$.  
%
  In the proof of Lemma~\ref{letri1} we pointed out that the following
  substitution is a bijection of Lawson quadrilaterals:
  \ben \label{glsubs}
    (l,r,s,t;\beta)\mapsto (l,r,\pi-s,\pi-t;\pi-\beta).
  \een
  When we take the quadrilaterals $\G_1$, this substitution gives 
  quadrilaterals $\G_2$ parameterized by~$\Z_2$ with $\pi/2<\beta<\pi$ and 
  unique lengths $\pi/2< s,t<\pi$.

  \noindent
  \ii\
  We now discuss the case $b=\alpha/2$ or $\beta=\pi/2$.  On $\Z^-_0$ we
  have $R(\alpha)<\pi/4$.  From~\eqref{gltrig} follows 
  $\cos^2\!2l=\cos^2\!2R(\alpha)$,
  so that $l=R(\alpha)$ by~\eqref{gllrra}.  This satisfies~\eqref{glbalr}
  too, thus the Lawson quadrilaterals are of the form
  $\Gamma(R(\alpha),t,R(\alpha),s;\pi/2)$.  Lemma~\ref{letri2}\ii\ gives the
  one-parameter family $\Gamma(R(\alpha),\pi,R(\alpha),\pi;\pi/2)$
  defining~$\G_0^-$.  By Lemma~\ref{leuniq} these are all such Lawson 
  quadrilaterals with $0<s,t\le \pi$.

  On $\Z^+_0$ we have $R(\alpha)>\pi/4$ and $l=\pi/2-R(\alpha)$
  satisfies~\eqref{glbalr}.  Now Lemma~\ref{letri2}\ia\ gives the
  quadrilaterals $\Gamma(\pi/2-R(\alpha),\pi/2,R(\alpha),\pi/2;\pi/2)$
  which are again unique.  These give the family~$\G_0^+$.

  Finally, at the point $\sigma$ with $R(\alpha)=l=\pi/4$ there is an entire 
  one-parameter family~$\G_0^{\sigma}$ of unique quadrilaterals 
  $\Gamma(\pi/4,s,\pi/4,s;\pi/2)$ with $0<s<\pi$ given by 
  Lemma~\ref{letri2}\iii. This proves \ii.

  We have to compare the quadrilaterals of \ii\ with the one-sided limits
  taken in $\Z_1$ and~$\Z_2$.  The points of $\Z_1$ converging
  to $\Z^-_0$ are of the form $(\alpha,r,b_1(\alpha,r))$ with
  $\alpha<\arccos(2/3)$.  Thus for $r\nearrow R(\alpha)$ we find $t\searrow
  0$ by~\eqref{glacht} and $s\searrow 0$ by~\eqref{glnull}.  On the other
  hand the limits from $\Z_2$ (given by the $b_2$ branch through the
  substitution~\eqref{glsubs} from the $b_1$ values) have $s,t \nearrow \pi$. 

  However, for the limits converging to $\Z^+_0$ we obtain continuity: 
  the limiting quadrilaterals with $(\alpha,r,b_1)$ for 
  $r\nearrow R(\alpha)>\pi/4$ have 
  $t\nearrow\pi/2$ by \eqref{glacht}, and $s\nearrow \pi/2$ by \eqref{glnull}.
  The same lengths are obtained for the similar limit with 
  $(\alpha,r,b_2)$ via~\eqref{glsubs}.
\end{proof}
  
%---------------------------------------------------------------------------
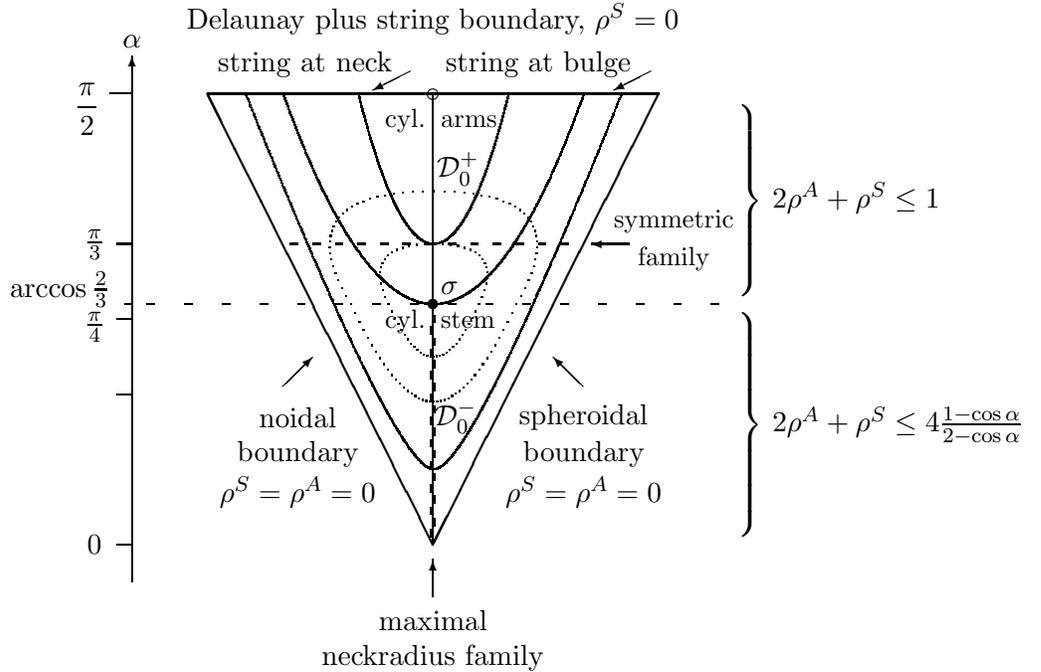
\begin{figure}[t]
  \small
  % \hspace*{-8mm}
    \unitlength=1mm
    \begin{picture}(120,100)(-20,5)
    % range 14 < x <  , 10 < y < 108
      %\put(15.00,20.00){\vector(1,0){70}}
      \put(10,87){\makebox(0,0)[cc]{$\alpha$}}
      \put(10.00,80.00){\line(-1,0){3.00}}
      \put(10.00,60.00){\line(-1,0){2.00}}
      \put(5,60.00){\makebox(0,0)[cc]{$\frac{\pi}3$}}
      \put(10.00,52.00){\line(-1,0){1}}
      \put(7,54){\makebox(0,0)[rc]{$\arccos\frac 2 3 $}}
      %\put(7,54){\makebox(0,0)[rc]{$\scriptsize\arccos\frac 2 3 $}}
      \put(10.00,50.00){\line(-1,0){2.00}}
      \put(5,49.00){\makebox(0,0)[cc]{$\frac{\pi}4$}}
      \put(10.00,40.00){\line(-1,0){2.00}}
      \put(10.00,20){\line(-1,0){2.00}}
      \put(10.00,15.00){\vector(0,1){70}}
      %\put(80.00,20.00){\line(0,-1){3.00}}
      %\put(60.00,20.00){\line(0,-1){2.00}}
      %\put(50,20.00){\line(0,-1){3}}
      %\put(40.00,20.00){\line(0,-1){2.00}}
      %\put(50,5){\makebox(0,0)[cc]{$\pi/2$}}
      \put(50,13){\vector(0,1){5}}
      \put(50,10){\makebox(0,0)[cc]{maximal}}
      \put(50,5){\makebox(0,0)[cc]{neckradius family}}
      %\put(80,15){\makebox(0,0)[cc]{$\pi$}}
      \put(4,78){\makebox(0,0)[cc]{$\displaystyle\frac{\pi}2$}}
      \put(5,20){\makebox(0,0)[cc]{$0$}}
      %\put(145,16.00){\makebox(0,0)[lc]{$\beta=\frac{\pi}2\!-\!\alpha\!+\!2b$}}
      %\put(87,20){\makebox(0,0)[cc]{$\beta$}}
      \qbezier(30,80)(40,52.00)(50.00,52.00)       % left,arms = 1/4
      \qbezier(70,80)(60,52)(50.00,52.00)             % right,arms = 1/4
      \qbezier(40,80)(45,60)(50,60)                    % left,arms = 1/3
      \qbezier(50,60)(55,60)(60,80)                    % right,arms = 1/3
      \qbezier(25,80)(45,30)(50,30)                    % left,arms = small
      \qbezier(75,80)(55,30)(50,30)                    % right,arms = small
      \qbezier[16](44,52)(40,60)(50,60)      % circle about sigma, up.left part
      \qbezier[16](56,52)(60,60)(50,60)      % circle about sigma, up.right part
      \qbezier[14](44,52)(47,45)(50,45)      % circle about sigma, l.left part
      \qbezier[14](56,52)(53,45)(50,45)      % circle about sigma, l.right part
      \qbezier[27](39,52)(30,67)(50,67)      % circle about sigma, up.left part
      \qbezier[27](61,52)(70,67)(50,67)      % circle about sigma, up.right part
      \qbezier[21](39,52)(45,39)(50,39)      % circle about sigma, l.left part
      \qbezier[21](61,52)(55,39)(50,39)      % circle about sigma, l.right part
      %\bezier{136}(35,80)(38,69)(55,54)                     % left,arms = 1/6
      %\bezier{124}(55,54)(71,40)(80,40)
      %\bezier{124}(80,40)(89,40)(105,54)                    % right,arms = 1/3
      %\bezier{136}(105,54)(122,69)(125,80)
      \put(50,52.00){\circle*{1.5}}
      \put(50,20.00){\line(0,1){60.00}}
      \multiput(31,60)(3,0){13}{\line(1,0){1}}
      \multiput(12,52)(5,0){16}{\line(1,0){1}}
    % \put(50,100){\makebox(0,0)[cc]{maximal}}
    % \put(50,95){\makebox(0,0)[cc]{family}}
      %\put(80,100){\makebox(0,0)[cc]{$r=R(\alpha)$,}}
      %\put(80,95){\makebox(0,0)[cc]{$b=\alpha/2$}}
    % \put(50,91){\vector(0,-1){7.00}}
    % \put(50,91){\vector(0,-1){5.00}}
      \put(70,37){\makebox(0,0)[cc]{spheroidal}}
      \put(70,32){\makebox(0,0)[cc]{boundary}}
      \put(70,27){\makebox(0,0)[cc]{$\rho^S=\rho^A=0$}}
      \put(70,41){\vector(-1,1){4.00}}
      \put(32,37){\makebox(0,0)[cc]{noidal}}
      \put(32,32){\makebox(0,0)[cc]{boundary}}
      \put(32,27){\makebox(0,0)[cc]{$\rho^S=\rho^A=0$}}
      \put(30,41){\vector(1,1){4.00}}
      \put(50,80.00){\circle{1.5}}
      % \put(49,83){\makebox(0,0)[rc]{cylindrical arms}}
   %  \put(51,83){\makebox(0,0)[lc]{$\rho^A=1/2$, $\rho^S=0$}}
      %\put(84,84){\vector(-1,-1){3}}
      \put(50,90){\makebox(0,0)[cc]{Delaunay plus string boundary, $\rho^S=0$}}
      \put(33,84){\makebox(0,0)[cc]{string at neck}}
      \put(47,83){\vector(-2,-1){5.00}}
      \put(64,84){\makebox(0,0)[cc]{string at bulge}}
      \put(79,83){\vector(-2,-1){5.00}}
      \put(49,78){\makebox(0,0)[rt]{\footnotesize cyl.}}
      \put(51,77){\makebox(0,0)[lt]{\footnotesize arms}}
      % \put(15.00,84.50){\makebox(0,0)[rb]{$\rho^A\in(0,1/2)$, $\rho^S=0$}}
    % \put(14,83){\vector(2,-1){5.00}}
      % \put(85,84.50){\makebox(0,0)[lb]{$\rho^A\in(0,1/2)$, $\rho^S=0$}}
    % \put(88,84){\vector(-2,-1){6.00}}
      \put(49,51.3){\makebox(0,0)[rt]{\footnotesize cyl.}}
      \put(51,51){\makebox(0,0)[lt]{\footnotesize stem}}
      \put(51,55){\makebox(0,0)[lt]{$\sigma$}}
      \put(76,60.00){\vector(-1,0){5.00}}
      \put(82,63){\makebox(0,0)[cc]{\footnotesize symmetric}}
      \put(82,58){\makebox(0,0)[cc]{\footnotesize family}} 
      \multiput(49.7,21)(0,3){10}{\line(0,1){1}}                  % for cut 2
      \multiput(50.2,22.5)(0,3){10}{\line(0,1){1}}                  % for cut 1
      \put(50.2,37){\makebox(0,0)[lc]{$\Z_0^-$}} 
      \put(50.5,70){\makebox(0,0)[lc]{$\Z_0^+$}} 
      \put(90,66){\makebox(0,0)[lc]{$\left.\rule{0mm}{14mm}\right\}$
	$2\rho^A+\rho^S  \le 1$}}
      \put(90,36){\makebox(0,0)[lc]{$\left.\rule{0mm}{16mm}\right\}$
	$2\rho^A+\rho^S \le 4 \frac{1-\cos\alpha}{2-\cos\alpha}$}}
      \thicklines
      \put(50,20.00){\line(-1,2){30.00}}
      \put(20.00,80.00){\line(1,0){60.00}}
      \put(80,80.00){\line(-1,-2){30.00}}
    \end{picture}
  \caption{  \label{fimois}
     The moduli space of the isosceles surfaces represented schematically
     as the interior of a triangle.  We indicate some curves of constant
     neckradius with solid lines for the arms, and dotted lines for the stem.
     Existence is only proved for the symmetric family~\cite{kgb} and in an
     infinite set in a neighborhood of the spheroidal boundary~\cite{kap}.
  }
\end{figure}
%---------------------------------------------------------------------------

%%%%%%%%%%%%%%%%%%%%%%%%%%%%%%%%%%%%%%
\subsection{All isosceles boundary contours}

The truncated pentagons obtained by applying Lemma~\ref{leistr}  
to the family of quadrilaterals constructed in Lemma~\ref{lemain} 
have an arc with (formal) length $s-b$, which may be negative.  
This is avoided by taking~$s$ sufficiently large.  We extend 
the arcs with lengths $b$, $t$, $s-b$ to infinity, and obtain
the boundary contours for isosceles surfaces directly from the quadrilaterals.
The infinite contours no longer have a discontinuity along~$\Z_0^-$.
%kgb

\begin{theorem} \label{thisum}
  An isosceles surface of angle $0<\alpha<\pi/2$ and with
  neckradii $(\rho^S,\rho^A)$ satisfies
  \ben
    &\rho^A \le \rho^A_{\max}(\alpha) := \frac{1-\cos\alpha}{2-\cos\alpha},
	  & \label{glrhoA}\\
    &\rho^S \le \rho^S_{\max}(\alpha) 
	 := \min\left\{  \frac{\cos\alpha}{2-\cos\alpha},
	     1-\frac{\cos\alpha}{2-\cos\alpha} \right\} &. \label{glrhoS}
  \een
  Furthermore, each fundamental domain of an isosceles surface is
  associated to a minimal domain in~$\s^3$ with boundary contained in some
  family~${\mathcal G}$, which is homeomorphic to a two-dimensional open
  disk.  On~$\partial \G$ the contours degenerate with $\rho^S=0$.  Each
  pair of neckradii satisfying the strict inequalities in~\eqref{glrhoA}
  and~\eqref{glrhoS} arises for two different boundary contours in~$\G$, and
  equality is satisfied with one contour.
%kgb
\end{theorem}
In particular the neckradius sum over all three ends satisfies
\ben \label{glnsum}
    2\rho^A+\rho^S  \le \min\left( 1,
    %\rho_1+\rho_2+\rho_3 & \le \min\left( 1,
    % \rho_1 + \rho_2 + \rho_3 & \le &
    4\,{\displaystyle \frac {1-\cos\alpha}{2-\cos\alpha}}\right) .
\een
\begin{proof}
  Extending $s$ and $t$ of the quadrilaterals in 
  ${\mathcal G}_1 \cup {\mathcal G}_2 \cup \G_0^{\pm}\cup\G_0^\sigma$ 
  by~$\pi$ and applying Lemma~\ref{leistr} yields
  non-degenerate pentagons.  We extend the arcs of the pentagon 
  with length $t$, $b$, and~$s-b$ to infinity.  This extension leads
  to the same boundary contour for all contours in $\G_0^\sigma$, and thus 
  all isosceles contours arising from Lemma~\ref{lemain} are 
  parameterized by~$\Z$.  Furthermore, they are continuous
  in $\Z$ since $l,r$ are continuous by Lemma~\ref{lemain}.
  These are all isosceles boundary contours, since truncation
  of a given isosceles contour leads to a pentagon for isosceles
  surfaces for which we can assume $0<s,t\le \pi$ if we move
  the perpendiculars; for such quadrilaterals, uniqueness is proven in
  Lemma~\ref{lemain}.
  
  Finally, on $\partial \Z$  we have $l\equiv 0$ so that $\rho^S=0$.  
  Thus~$\mathcal G$ is a maximal family.
\end{proof}

Again we expect a continuous family of isosceles \cmc surfaces to be in
a 1-1 relation to our boundary contours.  Existence is known in two
cases:
\begin{itemize}
  \item For an infinite set accumulating in the spheroidal boundary arc
    with $\rho^S=\rho^A=0$ (see~\ref{sssbdy} below) 
    \cite{kap}, and
  \item for the one-parameter subfamily of dihedrally symmetric
    triunduloids with $\alpha=\pi/3$ by~\cite{kgb}.
\end{itemize}

%%%%%%%%%%%%%%%%%%%%%%%%%%%%%%%%%%
\subsection{The geometry of the isosceles family}

On the assumption that our family of contours is bijective to a family of \cmc
surfaces, we discuss the geometrical properties of this family.
Using a numerical existence scheme, we found such a two parameter family
of isosceles surfaces, and give images in the paper~\cite{gbp}.

%%%%%%%%%%%%%%%%%%%%%
\subsubsection{The moduli space boundary} \label{sssbdy}

There are three geometrically distinguished arcs on the boundary of
isosceles contours, i.e.~the contours parameterized with~$\partial\Z$.
All contours on the boundary have vanishing neckradius of the stem end.

On the \emph{spheroidal boundary}, parameterized with 
$(\alpha,0,\alpha)\in\partial\Z$ we have $\rho^S=\rho^A=0$, so that the
limiting configuration consists of strings of spheres with isosceles
symmetry.  The limiting lengths (mod $\pi$) of the curvature arcs are
$t=\pi$, $b=\alpha$, and $s-b=\pi-\alpha$; considered as real numbers they
describe the lengths on the central junction sphere with punctures at the
limit points of the necks.  Kapouleas' method~\cite{kap} gives isosceles
surfaces, which accumulate at the spheroidal boundary arc.

On the \emph{noidal boundary}, parameterized with $(\alpha,0,0)$,
the ends are also spherical ($\rho^S=\rho^A=0$) but the contours 
have $b=t=0$.  This is the expected value for trinoidal junctions,
i.e.\ the three strings of spheres are attached to a point.  A blow-up
in the center of the \cmc surfaces such that the necks are scaled to constant
radius gives in the limit a minimal trinoid with isosceles symmetry.
The existence of
$k$-unduloids close to this noidal boundary component might follow from 
analogous work in the constant scalar curvature setting~\cite{mp}.

The noidal and spheroidal boundary agree in a limiting contour
with $\alpha=0$.  In this degenerate case the two arms coincide and
form a double string of spheres.

Finally there is an arc on~$\partial\Z$ with $\alpha=\pi/2$ which we call 
\emph{Delaunay plus string}.
Here only $\rho^S=0$, but $\rho^A\in(0,\pi/2]$.  For the half 
with $b<\alpha/2$, there is a string of spheres attached to a Delaunay
neck.  When $\rho^A\to 0$ the Delaunay surface itself tends to a string
of spheres, and the limit agrees with the noidal boundary case for
$\alpha=\pi/2$.  At the other endpoint, $\rho^A\to\pi/2$ the Delaunay
surface is a cylinder.  This leads over to the other half of the
family with $b>\alpha/2$ which is a similar family with the string of spheres 
attached at a Delaunay bubble.  Its limiting case is spheroidal.

%%%%%%%%%%%%%%%%%%%%%%%
\subsubsection{A triunduloid with one cylindrical end}

The contour at $\sigma\in\Z$ has  arm neckradius
$\rho^A=1/4$ and the stem has the cylindrical neckradius $\rho^S=1/2$.
The arms enclose an angle of $2\arccos(2/3)\approx 96.4\gd$.  This is the
only contour in the family with a cylindrical neckradius.  The arms can
attain the cylindrical neckradius only on the moduli space boundary.

%%%%%%%%%%%%%%%%%%%%%%%%%
\subsubsection{The maximal neckradius family}

This one-parameter subfamily of contours in $\mathcal G$ is parameterized by
$\Z_0$ and realizes the maximal neckradii for a given~$\alpha$.  It includes
the contour with cylindrical stem.  The right hand side of~\eqref{glrhoS}
gives the minimal and maximal asymptotic radius of the stem;  on $\Z^-_0$
the first radius is minimal, while on $\Z_0^+$ the second one is.  At the
cylindrical contour the minimum flips and the consequence is the
non-differentiability of the maximal neckradius sum~\eqref{glnsum}
at~$\sigma$:  the neckradius sum is~1 for $\arccos(2/3)\le\alpha<\pi/2$,
but smaller than~1 for $0<\alpha<\arccos(2/3)$.  If instead of the 
neckradius we replace~$\rho^S$ by the asymptotic bulgeradius
$\rho_0^S=\cos\alpha/(2-\cos\alpha)$ then $2\rho^A+\rho_0^S=1$ on the
maximal family.

It is interesting to observe that the asymptotic position of the minimum
neckradius shifts by~$\pi/2$ at $\sigma$, so that the selection of an
asymptotic planar curvature circle of radius~$\rho^S$ can not be analytic.
Indeed, our values $s\bmod\pi$, $t\bmod\pi$ given by Lemma~\ref{lemain} are
discontinuous on the maximal neckradius family at~$\sigma$, when they jump
by $\pi/2$, half a translational period of the end.

%%%%%%%%%%%%%%%%%%%%%
\subsubsection{Marked bubbles} \label{sssmar}

It may appear that selecting a family without requiring~\eqref{gllrra} 
would avoid the discontinuity of $t$ and~$s$ in our description, so that  
a perpendicular of the boundary rays could be marked in a continuous way
over the family.  However, we want to show that a discontinuity must be
present around~$\sigma$.

Let $\gamma\subset\Z-\{\sigma\}$ be a loop  winding once about the
point~$\sigma$.  The isosceles boundary contours are continuous
on~$\gamma$.  Nevertheless, according to Lemma~\ref{lemain}, the
corresponding path of quadrilaterals~$\Gamma_L(\gamma)$ has a discontinuity
of~$s,t$ of~$\pm\pi$, when $\Z_0^-$ is crossed.  Thus the asymptotic
perpendicular moves by one full period along the stem.  This means that on
the closed loop~$\gamma$ a bubble is created or deleted on the isosceles
surfaces.  Note also that the discontinuity of $t,s$ can be located on a
curve joining any point of~$\partial\Z$ to~$\sigma$.

As a consequence, a closed curve winding $n$~times about~$\sigma$ generates
or deletes $n$~bubbles on the stem.  For all isosceles surfaces, except for
the one with cylindrical stem, we can mark a bubble on the stem.  The
corresponding ``marked'' moduli space is then the universal covering of the
annulus $\Z-\{\sigma\}$.

%%%%%%%%%%%%%%%%%%%%%%%%%%%%%%%%%%%%%%%%%%%%%%%%%%%%%%%%%%%%%%%%
%\note{How intrinsic is the spherical geometry?}
%
% The present paper gives evidence that the impact of spherical geometry
% on the behavior of \cmc surfaces is rather strong.  Surprisingly,
% the geometrically defined weights seem less important than necksizes:
% for instance, if we express~\eqref{glrsum} in terms of weights
% we obtain
  % a = r_min/r_max = r_min/(1-r_min) = 1 / [ 1/r_min - 1 ]
  % ==> a_1 + a_2 = some computation
  %                     .5 - 2 r + 2 r^2
  %                 \le ----------------   with r = r_min,1
  %                     .5 + .5 r - r^2
% $$ \sqrt {\frac 1 4 - \frac{w_1}{2\pi}}\,
  % +  \sqrt {\frac 1 4 - \frac{w_2}{2\pi}} \ge \frac 1 2 .
% $$
% Furthermore the special property of a cylindrical end is that its
% associated boundary curves have not a discrete but a continuous set of 
% perpendiculars, see the distinction made in Lemma~\ref{lequpa}.
% For finite length segments, the isoperimetric inequality gives that 
% this property 
% occurs at weights larger than a cylinder 
% (for an example of a finite length cylindrical piece see \cite{gbp}).
% We would also remark that the branching of rhombic surfaces with 
% necksize~1/4 (see~\cite{gbp})is an immediate consequence of the
% branching of right-angled Lawson quadrilaterals given in Lemma~\ref{letri2}.

%%%%%%%%%%%%%%%%%%%%%%%%%%%%%%%%%%%%%%%%%%%%%%%%%%%%%%%%%%%%%%%%%%%%%%%%%%

\end{document}